\documentclass[12pt]{elsarticle}

\usepackage{graphicx}

\usepackage{subcaption}
\usepackage{amsmath}
\usepackage{enumitem}
\numberwithin{equation}{section}
\usepackage{algorithm,algorithmic}
\usepackage{bm}
\usepackage{amssymb}
\usepackage{geometry}
\begin{document}

\begin{frontmatter}



\title{R-PINN: Recovery-type a-posteriori estimator enhanced adaptive PINN} 
\author[label1]{Rongxin Lu} 
\author[label1,label2]{Jiwei Jia\corref{cor1}}
\cortext[cor1]{Corresponding author.}
\ead{jiajiwei@jlu.edu.cn}
\author[label3]{Young Ju Lee}
\author[label1]{Zheng Lu}
\author[label4,label5]{Chen-Song Zhang}

\address[label1]{Department of Computational Mathematics, School of Mathematics, Jilin University, Changchun, China}
\address[label2]{AI for Science and Engineering Center, Shenzhen Loop Area Institute, Shenzhen 518048, China}
\address[label3]{Department of Mathematics, Texas State University, San Marcos, TX, USA}
\address[label4]{SKLMS and NCMIS, Academy of Mathematics and Systems Science, Beijing 100190, China} 
\address[label5]{School of Mathematical Sciences, University of Chinese Academy of Sciences, Beijing 100049,
China}

\begin{abstract}
	In recent years, with the advancements in machine learning and neural networks, algorithms using physics-informed neural networks (PINNs) to solve PDEs have gained widespread applications. While these algorithms are well-suited for a wide range of equations, they often exhibit a suboptimal performance when applied to equations with large local gradients, resulting in substantially localized errors. To address this issue, this paper proposes an adaptive PINN algorithm designed to improve accuracy in such cases. The core idea of the algorithm is to adaptively adjust the distribution of collocation points based on the recovery-type a-posteriori error of the current numerical solution, enabling a better approximation of the true solution. This approach is inspired by the adaptive finite element method. By combining the recovery-type a-posteriori estimator, a gradient-recovery estimator commonly used in the adaptive finite element method (FEM), with PINNs, we introduce the recovery-type a-posteriori estimator enhanced adaptive PINN (R-PINN) and compare its performance with a typical adaptive sampling PINN, failure-informed PINN (FI-PINN), and a typical adaptive weighting PINN, residual-based attention in PINN (RBA-PINN) as a baseline. Our results demonstrate that R-PINN achieves faster convergence with fewer adaptively distributed points and outperforms the other two PINNs in the cases with regions of large errors.
	
\end{abstract}




\end{frontmatter}
\section{Introduction} 

The physics-informed neural networks (PINNs) algorithm is a data-driven approach for solving partial differential equations (PDEs) based on deep learning frameworks \cite{Raissi2019}. Due to its distinct methodology, which diverges significantly from traditional numerical techniques, and by leveraging the superior capabilities of neural networks in handling large-scale data and optimization problems, PINNs have garnered considerable attention ever since their introduction and have become important tools for solving PDEs. These methods utilize automatic differentiation to compute the residual of an equation at a set of collocation points within the domain of interest, using it as an objective loss function that is optimized through backpropagation, thereby approximating the solution of the PDE. Compared to classical PDE numerical methods such as the finite element method and finite difference method, one significant advantage of PINNs is that they are mesh-free, making them easier to implement and particularly well-suited for high-dimensional problems and irregular domains. Currently, PINNs have been applied to solve a range of PDE problems, including fluid mechanics, which has been reviewed comprehensively by Cai et al.\  \cite{cai2021physics}, as well as  fields such as optics, metamaterials and others \cite{raissi2020hidden, chen2020physics, pang2019fpinns, zhang2019quantifying, guo2022normalizing}.

Despite the notable advantages demonstrated by PINNs in PDEs, there remains a substantial room for improvement. While PINNs have exhibited high-precision solutions across a variety of scenarios, their applicability is still limited. Notably, significant errors can arise when addressing equations such as the Allen-Cahn equation or the Burgers' equation. Enhancements to the PINNs methodology can be pursued from multiple aspects. For instance, since the total loss function in PINNs is a weighted sum of errors from both the domain and boundary conditions, adjusting the weights of these terms may improve performance across different problems \cite{ANAGNOSTOPOULOS2024116805,wang2021understanding, lu2021physics, wang2022and, gu2021selectnet}. A representative example is the residual-based attention in PINN \cite{ANAGNOSTOPOULOS2024116805}, which  selects weighting multipliers for the loss function and has been shown to significantly enhance the performance of PINNs. In cases involving extensive solution domains,  decomposition techniques can be employed to localize the training process, thereby enhancing optimization outcomes \cite{meng2020ppinn, shukla2021parallel, jagtap2020extended}. Regarding the choice of loss functions, although the mean squared error (MSE) is commonly used, alternative forms of loss functions may give better training results in certain circumstances \cite{gu2021selectnet,psaros2022meta, mcclenny2020self}. Furthermore, strategies, such as batch training and the design of specialized network architectures to deal with Dirichlet or periodic boundary conditions, offer a promising direction to enhance the performance of PINNs \cite{lu2021physics, lagari2020systematic, dong2021method}.

In this paper, we focus on the limitations of PINNs in handling solutions with sharp local variations. This challenge, arising from singularities or high-gradient regions, is typically addressed by introducing adaptive sampling. However, existing sampling methods primarily rely on the PDE loss at local collocation points, which can be regarded as a residual-type estimator. Nevertheless, alternative efficient estimators remain unexplored for adaptive PINNs. Therefore, we propose the recovery-type a-posteriori estimator enhanced adaptive PINN (R-PINN), which integrates the recovery-type a-posteriori estimator from classical numerical methods into PINNs. Furthermore, we introduce an adaptive sampling method, the recovery-type estimator adaptive distribution (RecAD), which performs adaptive refinement based primarily on the recovery-type error estimator.

The sampling methods for collocation points in PINNs are generally categorized into uniform and non-uniform sampling. In typical uniform sampling, points are either randomly distributed or placed at equispaced intervals, remaining fixed throughout the training process. However, this static sampling strategy has notable limitations, as it fails to capture the model's fine local features effectively, which can lead to convergence issues in complex problems settings. To address these limitations, Lu et al. first proposed an adaptive non-uniform sampling method \cite{lu2021deepxde}. Since then, various adaptive sampling methods for collocation points have been developed \cite{wu2023comprehensive,gao2023failure,nabian2021efficient,gao2023active,zeng2022adaptive,hanna2022residual}. Among these, the failure-informed PINN (FI-PINN) by Gao et al. \cite{gao2023failure} is particularly notable, utilizing a truncated Gaussian model for adaptive sampling. Another significant contribution is by Wu et al. \cite{wu2023comprehensive}, who reviewed and compared a range of sampling techniques. Their work covers both uniform sampling methods, such as equispaced, random, and low-discrepancy sequences, and non-uniform methods, including residual-based adaptive refinement (RAR) and residual-based adaptive distribution (RAD), offering valuable guidance on selecting appropriate sampling strategies.

While there are various adaptive PINNs methods, they generally share a common feature: sampling the uniform points as background collocation points and then adding non-uniform points as the adaptively distributed points based on an error estimator. Most of these methods utilize the residual-type a-posteriori estimator, which calculates the error as the loss of the neural network function. For adaptive finite element methods, two primary types of estimator are widely recognized: residual-type \cite{babuvska1978error,babuvska1978posteriori} and recovery-type \cite{Zienkiewicz1987,zienkiewicz1992superconvergent1,zienkiewicz1992superconvergent2}, with the latter based on gradient recovery. Although extensive research has been conducted on residual-type adaptive PINNs \cite{wu2023comprehensive,gao2023failure}, the use of gradient recovery-based estimators remains unexplored. Therefore, we introduce an adaptive PINN enhanced by the recovery-type a-posteriori estimator.

The rest of this paper is organized as follows. Section 2 provides preliminary background, introducing the fundamental algorithm of physics-informed neural networks (PINNs) and existing sampling strategies for both uniform and non-uniform purposes. Section 3 presents the detailed algorithm of the proposed recovery-type estimator-enhanced PINN (R-PINN). Section 4 presents our results using R-PINN on four representative singular problems, illustrating its performance and providing a comparative analysis. Finally, Section 5 summarizes our conclusions and outlines potential directions for future research.

\section{Preliminaries}
In this section, we will present the requisite preliminary knowledge. We will begin with an overview of PINNs algorithm in solving both forward and inverse PDEs problems. Following this, we will discuss some strategies for collocation points sampling employed in the implementation of PINNs, including both uniform and non-uniform sampling methodologies.

\subsection{PINNs for solving forward and inverse PDEs}
Analogous to the application of basis functions for approximating true solutions in traditional numerical methods for PDEs, the physics-informed neural networks (PINNs) algorithm utilizes deep neural networks to approximate the exact solutions of equations, primarily employing fully connected feed-forward neural networks. We denote a fully connected feed-forward neural network with $L$ layers as $\mathcal{N}^L:\mathbb{R}^{N_0}\rightarrow\mathbb{R}^{N_k}$, where the $k$-th layer comprises $N_k$ neurons. The weight matrix and bias vector of the $k$-th layer (for $1\leq k\leq L$) can be represented as $\bm{W}^k\in\mathbb{R}^{N_k\times N_{k-1}}$ and $\bm{b}^k\in\mathbb{R}^{N_k}$, respectively. The input vector is denoted as $\bm{x}\in\mathbb{R}^{N_0}$, and the output vector of the $k$-th layer is expressed as 
$\mathcal{N}^k(\bm{x})$, especially, $\mathcal{N}^0(\bm{x})=\bm{x}$. The activation function $\Phi$, applied during the feed-forward transition from the $(k-1)$-th layer to the $k$-th layer (for $1\leq k\leq L$), introduces non-linearity into the network, thereby enhancing its learning capacity. Consequently, a feed-forward neural network can be mathematically defined as
\begin{equation*}
\mathcal{N}^k(\bm{x})=\Phi(\bm{W}^k\mathcal{N}^{k-1}(\bm{x})+\bm{b}^k),\quad 1\leq k\leq L.
\end{equation*}
Letting $\bm{\Theta}=\{\bm{W}^k, \bm{b}^k\}\in\mathcal{V}$ be all parameters of the neural network, including weights and biases, where $\mathcal{V}$ denotes the parameter space, we can finally express the neural network as
\begin{equation*}
    \hat{u}_{\bm{\Theta}}(\bm{x})=\mathcal{N}^L(\bm{x};\bm{\Theta}),
\end{equation*}
where $\mathcal{N}^L(\bm{x};\bm{\Theta})$ emphasizes the dependence of the output layer $\mathcal{N}^L(\bm{x})$ on $\bm{\Theta}$. Initial parameters are typically generated according to a specified probability distribution.

In the context of solving PDEs with PINNs, we consider a PDE defined over the domain $\Omega\in\mathbb{R}^d$ with parameters $\bm{\lambda}$,
\begin{equation*}
    f(\bm{x}; u(\bm{x})) = f\left(\bm{x}; \frac{\partial u}{\partial x_1}, \cdots, \frac{\partial u}{\partial x_d}; \frac{\partial^2 u}{\partial x_1^2}, \cdots, \frac{\partial^2 u}{\partial x_1 \partial x_d}; \cdots; \bm{\lambda}\right) = 0, \quad \bm{x} = (x_1, x_2, \cdots, x_d) \in \Omega, 
\end{equation*}
subject to the boundary condition:
\begin{equation*} 
    \mathcal{B}(u, \bm{x}) = 0 \quad \text{on} \ \partial\Omega,
\end{equation*}
where $u(\bm{x})$ represents the value of the true solution at point $\bm{x}$. 

In solving the forward problem, where the parameter $\bm{\lambda}$ in the equation is known, PINNs approximate the true solution $u$ using a neural network $\hat{u}_{\bm{\Theta}}(\bm{x})$  \cite{Raissi2019}. During this process, the neural network is iteratively trained to update its parameters $\bm{\Theta}$ in order to minimize the loss function, defined as,
\begin{equation}\label{loss1}
\mathcal{L}(\mathcal{T};\bm{\Theta})=w_f\mathcal{L}_f(\mathcal{T}_f;\bm{\Theta})+w_b\mathcal{L}_b(\mathcal{T}_b;\bm{\Theta}),
\end{equation}
where
\begin{equation*}
\begin{split}
    \mathcal{L}_f(\mathcal{T}_f;\bm{\Theta})=\frac{1}{\lvert\mathcal{T}_f\rvert}\sum_{\bm{x}\in\mathcal{T}_f}\bigg| &f\left(\bm{x}; \frac{\partial\hat{u}_{\bm{\Theta}}}{\partial x_1}, \cdots, \frac{\partial\hat{u}_{\bm{\Theta}}}{\partial x_d}; \frac{\partial^2\hat{u}_{\bm{\Theta}}}{\partial x_1\partial x_1},\cdots,\frac{\partial^2\hat{u}_{\bm{\Theta}}}{\partial x_1\partial x_d}; \cdots; \bm{\lambda}\right)  \bigg|^2,\\
    &\mathcal{L}_b(\mathcal{T}_b;\bm{\Theta})=\frac{1}{\lvert\mathcal{T}_b\rvert}\sum_{\bm{x}\in\mathcal{T}_b}\lvert\mathcal{B}(\hat{u}_{\bm{\Theta}},\bm{x})\rvert^2,
\end{split}
\end{equation*}
and $w_f$ and $w_b$ are weighting factors. Here, $\mathcal{T}_f$ and $\mathcal{T}_b$ denote the sets of sampled points within the domain and on the boundary, respectively, which are commonly referred to as "collocation points". The combined set of collocation points is given by $\mathcal{T}=\mathcal{T}_f\cup\mathcal{T}_b$.

For certain problems, while the parameters $\bm{\lambda}$ in the governing equations remain unknown, the solutions are observable, thereby giving rise to inverse problems. To solve the inverse problem with PINNs, the loss function not only accounts for the deviations dictated by physical laws but also incorporates an additional term to capture the residual between observed solution $u(\bm{x})$ and predicted solution $\hat{u}(\bm{x})$ over a set of collocation points $\mathcal{T}_i$:
\begin{equation*}
    \mathcal{L}_i(\mathcal{T}_i; \bm{\Theta}, \bm{\lambda}) = \frac{1}{|\mathcal{T}_i|} \sum_{\bm{x} \in \mathcal{T}_i} |\hat{u}_{\bm{\Theta}}(\bm{x}) - u(\bm{x})|^2.
\end{equation*}
As a result, the total loss function for the inverse problem can be expressed as:
\begin{equation}\label{loss2}
    \mathcal{L}(\mathcal{T}; \bm{\Theta}, \bm{\lambda}) = w_f \mathcal{L}_f(\mathcal{T}_f; \bm{\Theta}, \bm{\lambda}) + w_b \mathcal{L}_b(\mathcal{T}_b; \bm{\Theta}, \bm{\lambda}) + w_i \mathcal{L}_i(\mathcal{T}_i; \bm{\Theta}, \bm{\lambda}),
\end{equation}
with the weight of the additional term $w_i$. Then the unknown parameters $\bm{\lambda}$ can be trained with the neural network parameters $\bm{\Theta}$ by the optimization process.

Note that in equations~(\ref{loss1}) and (\ref{loss2}), it is crucial to position the set of collocation points within the domain of interest. While simple sampling methods, such as uniform and random sampling, are available, they are generally less effective in capturing the solution's local features. Consequently, in certain cases, employing more advanced non-uniform sampling techniques becomes particularly important.

\subsection{Sampling methods for PINNs}\label{sampling-method}
\begin{figure}[htbp]
    \centering
    \includegraphics[height=0.2\textheight]{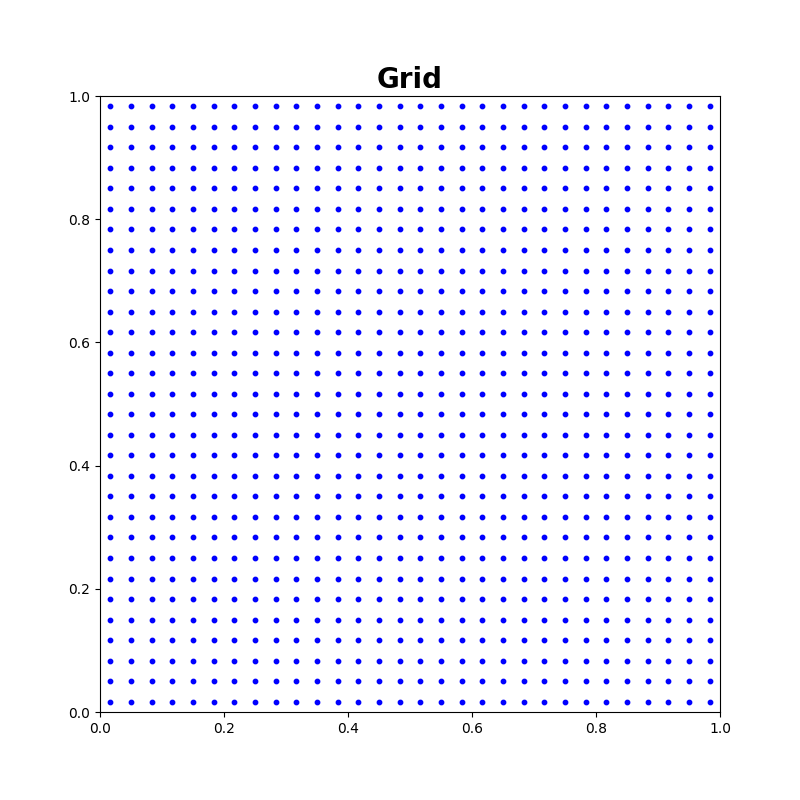}
    \includegraphics[height=0.2\textheight]{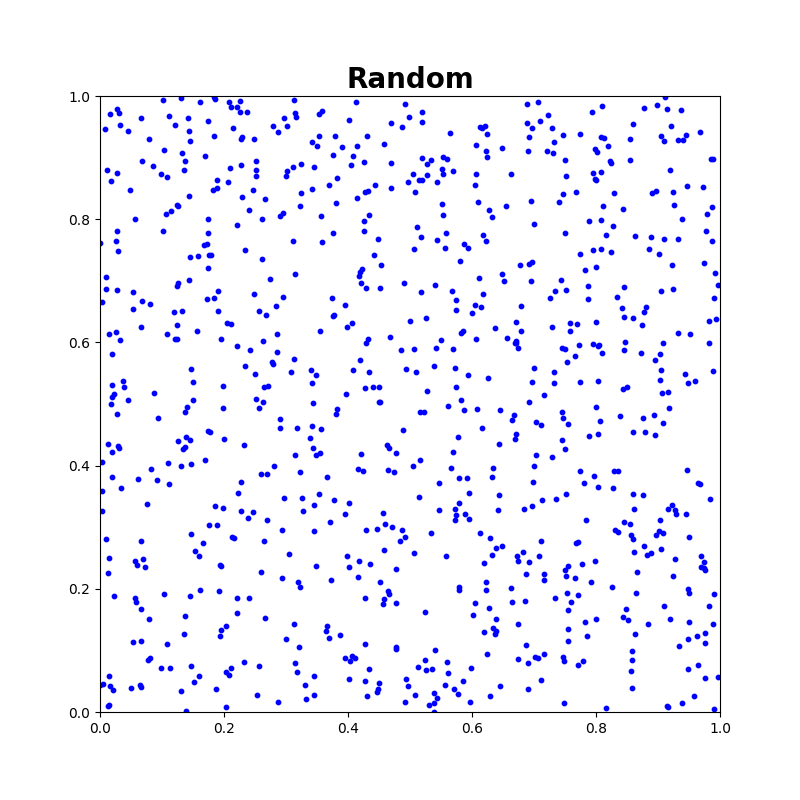}
    \includegraphics[height=0.2\textheight]{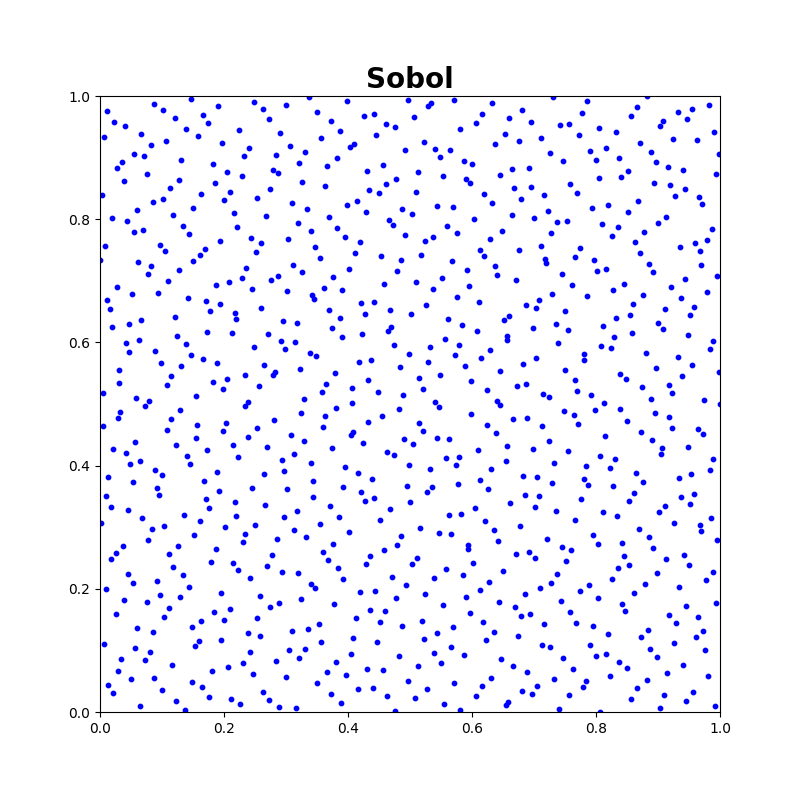}
    \caption{$900$ uniform collocation points generated on $[0, 1]^2$ with three different strategies.}
    \label{uniform-points}
\end{figure}

Typically, the sampling strategy in PINNs involves selecting a set of collocation points by sampling prior to training the neural network, with the set remaining fixed throughout the training process. 
Common strategies for uniformly sampling collocation points include equispaced grid sampling and uniformly random sampling. Several other quasi-uniform sampling techniques have also been proposed in the literature \cite{Raissi2019,pang2019fpinns,guo2022analysis}.
In  Wu et al.'s work \cite{wu2023comprehensive}, six uniform sampling methods were systematically reviewed, including equispaced uniform grid (Grid), uniformly random sampling (Random), Latin hypercube sampling (LHS) \cite{mckay2000comparison,stein1987large}, Halton sequence (Halton) \cite{halton1960efficiency}, Hammersley sequence (Hammersley) \cite{hammersleymonte}, and Sobol sequence (Sobol) \cite{sobol1967distribution}. Their comparative experiments demonstrated that the Sobol method generally provides superior performance. The Sobol sequence uses a base of two to form successively finer uniform partitions of the unit interval and then reorders the coordinates in each dimension. We show the difference over three different methods in Fig.~\ref{uniform-points}.

In many cases, accuracy of numerical solution can be good enough through uniform sampling strategies. However, for certain numerical problems characterized by large local gradients and highly distinctive solution features, such methods may fail to adequately capture the local variations. Several adaptive methods have been proposed to overcome these limitations in recent years, primarily through modifications to the sampling strategies for the collocation points in PINNs. These adaptive sampling approaches can be broadly classified into two categories:
\begin{itemize}
    \item The first category incrementally adds newly selected points based on an adaptive criterion, meaning that previously selected points are retained and the number of sampling points increases over time. Representative methods in this category include the residual-based adaptive refinement (RAR) \cite{lu2021deepxde}, which identifies regions with large residuals and adds new points periodically, and the self-adaptive importance sampling (SAIS) \cite{gao2023failure}, which constructs a probability density function (PDF) based on the residual to sample new points and adds them to the existing dataset.
    \item The second category, in contrast, maintains a fixed total number of sampling points by periodically redistributing them based on an updated sampling criterion. A typical example is the residual-based adaptive distribution (RAD) \cite{nabian2021efficient}, which uses a residual-informed PDF $p(\bm{x})$ to resample all collocation points. A recent study \cite{wu2023comprehensive} shows that methods of the second category often achieve better performance. Therefore, in this work, we adopt the strategy of shifting all adaptively distributed points periodically.
\end{itemize}

The adaptive sampling methods discussed above are all based on the residual-type error estimator. However, another practical alternative, the recovery-type error estimator, has not yet been explored in PINNs-based methods. By relying on reconstructed solution fields rather than direct residuals, recovery-based estimators provide a different perspective on error characterization, which can lead to more robust and informative sampling criteria. This motivates us to develop an adaptive PINN framework based on recovery-type error estimation, along with a corresponding sampling strategy.

\section{Recovery-type estimator enhanced PINN}
In this section, we present an adaptive PINN algorithm, which utilizes the recovery-type a-posteriori error estimator \cite{Zienkiewicz1987,zienkiewicz1992superconvergent1,zienkiewicz1992superconvergent2}, a widely used technique in adaptive finite element methods. The concept of adaptive finite element method, first introduced by I. Babuška et al. in 1978 \cite{babuvska1978error}, has since seen significant advancements and remains integral to practical scientific computations. This approach extends classical finite element methods by incorporating a-posteriori error estimation and automated mesh refinement to iteratively improve the accuracy of approximated solutions, which can achieve the target accuracy with reduced computational effort, offering distinct advantages over traditional finite element methods, particularly when addressing problems involving singular solutions.

The fundamental concept of the adaptive finite element method lies in the use of a-posteriori error estimator to assess the numerical solution error within each mesh element, thereby guiding the mesh refinement and coarsening process. This adaptive mesh strategy shares a notable similarity with adaptive PINNs, where the distribution of collocation points is adjusted based on the error distribution. In adaptive finite element methods, commonly employed a-posteriori error estimators can be categorized into two main types: the first is based on residual type \cite{babuvska1978error,bank1985}, where the error is estimated by computing the discrepancy between the solution and the residual of the governing equation; the second type, recovery type, proposed by Zienkiewicz and Zhu \cite{Zienkiewicz1987,zienkiewicz1992superconvergent1,zienkiewicz1992superconvergent2}, is a gradient reconstruction-based method that provides more accurate error estimates by reconstructing the solution gradients.

In existing adaptive PINNs algorithms, the PDE residual, computed using automatic differentiation, is typically used as the default a-posteriori error estimator, which corresponds to the residual-based error estimator in adaptive finite element methods. Our proposed algorithm seeks to introduce the recovery-type error estimator into the adaptive PINNs framework, with the aim of enriching and extending the theoretical foundation of adaptive PINNs, ultimately providing a more comprehensive approach to error estimation and control.

\subsection{Recovery-type a-posteriori estimator \label{estimator}}

In this section, we provide a brief overview of recovery-type a-posteriori error estimators in the context of finite element methods, along with a discussion of how they work in practice. A more detailed theoretical treatment can be found in \cite{xu2004analysis, verfurth2013posteriori}.

Let $\mathcal{T}_h$ be a quasi-uniform triangulation of the domain $\Omega$, and denote by $\mathcal{Z}_h$ the set of its vertices. For a given node $z \in \mathcal{Z}_h$, we define an element patch $\omega$ consisting of all elements that share the vertex $z$. We place the origin of a local coordinate system at $z$, and denote by $(x_j, y_j)$ the barycenter of each triangle $\tau_j \subset \omega$, where $j = 1, 2, \dots, m$.

As a model problem, we consider the weak form of a  non-self-adjoint problem: Find $u \in H^1(\Omega)$ such that
\begin{equation*}\label{bilinear}
B(u, v) = \int_\Omega \left[ (\mathcal{D} \nabla u + \mathbf{b} u) \cdot \nabla v + c u v \right]\, \mathrm{d}{\bm{x}} = f(v), \quad \forall v \in H^1(\Omega),
\end{equation*}
with a Dirichlet boundary condition $u|_{\partial\Omega}=g$,
where $\mathcal{D}$ is a symmetric positive-definite $2 \times 2$ matrix, $\mathbf{b}\in [L^\infty (\Omega)]^2$ is a vector-valued function, $f\in (H^1(\Omega))'$ is a linear functional, and $g\in H^{\frac{1}{2}}(\partial\Omega)$ is a prescribed function defined on the boundary. We assume that the coefficients are smooth, and that the bilinear form $B(\cdot,\cdot)$ is continuous and satisfies an inf-sup condition on $H^1(\Omega)$, ensuring the well-posedness of the problem.

The finite element approximation $u_h \in \mathcal{V}_h \subset H^1(\Omega)$ is defined by
\begin{equation*}\label{finite-prob}
B(u_h, v_h) = f(v_h), \quad \forall v_h \in \mathcal{V}_h.
\end{equation*}
To guarantee a unique discrete solution, we further assume that the inf-sup condition also holds on the discrete space $\mathcal{V}_h$.

To indicate the recovery-type a-posteriori error, we introduce a gradient recovery operator $G$ designed to provide an improved approximation of the true gradient $\nabla u$. Ideally, the recovered gradient $G u_h$ satisfies the following inequality:
\begin{equation*}
\lVert\nabla u - G u_h\rVert_{0,\Omega} \le \beta \lVert\nabla u - \nabla u_h\rVert_{0,\Omega}, \quad \text{with } 0 \le \beta < 1,
\end{equation*}
which leads to a two-sided bound:
\begin{equation*}
(1 - \beta)\lVert\nabla u - \nabla u_h\rVert_{0,\Omega} \le \lVert\nabla u_h - G u_h\rVert_{0,\Omega} \le (1 + \beta)\lVert\nabla u - \nabla u_h\rVert_{0,\Omega}.
\end{equation*}
This justifies the use of $\lVert\nabla u_h - G u_h\rVert_{0,\Omega}$ as a recovery-based a-posteriori error estimator—it is fully computable and reliably reflects the true error.

There are several standard approaches for constructing $G u_h$ at a vertex $z$:
\begin{enumerate}[label=(\alph*)]
\item Weighted averaging:
The recovered gradient is constructed as a weighted average of the gradients evaluated at the barycenters of the neighboring triangles:
\begin{equation*}\label{weight}
G u_h(z) = \sum_{j=1}^m \frac{|\tau_j|}{|\omega|} \nabla u_h(x_j, y_j).
\end{equation*}

\item Local $L^2$ projection:
For each component $l = 1, 2$, we seek a linear function $p_l \in P_1(\omega)$ satisfying
\begin{equation*}\label{L2-proj}
\int_\omega \left[p_l(x, y) - \partial_l u_h(x, y)\right] q(x, y)\, \mathrm{d}x\mathrm{d}y = 0 \quad \forall q \in P_1(\omega),
\end{equation*}
and define $G u_h(z) = (p_1(0, 0), p_2(0, 0))$.
\item Local least-squares fitting (Zienkiewicz-Zhu method \cite{zienkiewicz1992superconvergent1}):
This method solves a discrete version of the $L^2$ projection, by minimizing the difference at the barycenters:
\begin{equation*}\label{least-square}
\sum_{j=1}^m \left[p_l(x_j, y_j) - \partial_l u_h(x_j, y_j)\right] q(x_j, y_j) = 0, \quad \forall q \in P_1(\omega), \quad l = 1, 2.
\end{equation*}
Again, we define $G u_h(z) = (p_1(0, 0), p_2(0, 0))$.
\end{enumerate}

Subsequently, $Gu_h$ can be linearly interpolated from $Gu_h(z),\ z\in\mathcal{Z}_h$. Among these three approaches, the weighted averaging method is particularly popular due to its simplicity and low computational cost. The least-squares fitting can be viewed as a discrete analogue of the $L^2$ projection. Detailed analysis on the convergence and performance of these estimators can be found in \cite{zienkiewicz1992superconvergent1, zienkiewicz1992superconvergent2, xu2004analysis, verfurth2013posteriori}.

\subsection{Algorithm of the recovery-type estimator enhanced PINN}
In this section, we introduce the algorithm of recovery-type estimator enhanced adaptive PINN (R-PINN). The success of the adaptive PINNs algorithm hinges on two crucial components: the a-posteriori error estimation method and the adaptive sampling strategy. The a-posteriori error estimation predicts the discrepancy between the numerical and true solutions, while the adaptive sampling strategy adjusts the distribution of collocation points based on this error, allowing for more focused and efficient training on areas that need it the most.

We begin by discussing a-posteriori error estimation based on the recovery-type estimator. The core idea is to utilize the estimator introduced in Section~\ref{estimator}, particularly the weighted averaging scheme, to assist in evaluating the a-posteriori error. The procedure is as follows. After training the PINN, we obtain an approximated solution $\hat{u}_{\bm{\Theta}} = \mathcal{N}^L(\bm{x}; \bm{\Theta})$. Given a quasi-uniform triangulation $\mathcal{T}_h$ of the domain $\Omega$, we interpolate $\hat{u}_{\bm{\Theta}}$ at mesh vertices to obtain a piecewise linear approximation $\hat{u}_I$ such that
\begin{equation*}
    \hat{u}_I(z)=\hat{u}_{\bm{\Theta}}(z),\quad \forall z\in\mathcal{Z}_h,
\end{equation*}
where $\mathcal{Z}_h$ denotes the set of mesh vertices. The interpolated solution $\hat{u}_I$ serves as an intermediate quantity for both gradient recovery and subsequent a-posteriori error estimation. The recovered gradient at a node $z \in \mathcal{Z}_h$ is computed using a weighted averaging approach,
\begin{equation*}
    G\hat{u}_I(z)=\sum_{j=1}^m\frac{|\tau_j|}{|\omega|}\nabla \hat{u}_I(x_j,y_j),
\end{equation*}
where $(x_j, y_j)$ is the barycenter of each triangle $\tau_j\subset\omega$, and $\omega$ is an element patch consisting of all elements that share the vertex $z$. 
The recovered gradient field $G\hat{u}_I$ is then obtained by linear interpolation of the nodal values $G\hat{u}_I(z)$. Using the recovered gradient, the recovery-type a-posteriori error estimator is defined element-wise as
\begin{equation}\label{estimator-eq}
    \eta_K=\lVert\nabla\hat{u}_I-G\hat{u}_I\rVert_{0,K},\quad K\in\mathcal{T}_h.
\end{equation}
Here, $\eta_K$ provides a local indication of the error within element $K$, guiding the identification of regions where additional sampling points are required to enhance the solution accuracy.

Before presenting the sampling strategy, we outline the overall training strategy. The training process is divided into two phases: a pre-training phase and an adaptive training phase.

In the pre-training phase, the PINN is trained using only a fixed set of background collocation points, denoted by  $\{\bm{x}_b^{(i)}\}_{i=1}^{N_1}$. This phase aims to obtain an initial approximation of the solution. The pre-training is followed by an adaptive training phase, during which both background collocation points $\{\bm{x}_b^{(i)}\}_{i=1}^{N_1}$ and a set of adaptively distributed points $\{\bm{x}_a^{(i)}\}_{i=1}^{N_2}$ are employed. While the background collocation points remain fixed throughout the entire training process, the adaptively distributed points are dynamically updated in each iteration. Specifically, the adaptive phase comprises $M$ iterations. At the beginning of each iteration, the recovery-based a-posteriori error estimator $\eta_K$ is computed to evaluate the current approximation error over each element $K$. Based on these error indicators, new adaptively distributed points are selected in regions with relatively large errors. The PINN is then retrained using the union of background and updated adaptively distributed points for a fixed number of epochs before proceeding to the next iteration.

During the adaptive training process, the sampling strategy plays a key role in determining how the collocation points are selected and updated. The background collocation points, denoted by $\{\bm{x}_b^{(i)}\}_{i=1}^{N_1}$, are generated using a Sobol sequence as described in Section~\ref{sampling-method}, and remain fixed throughout the entire training procedure. In contrast, the adaptively distributed points $\{\bm{x}_a^{(i)}\}_{i=1}^{N_2}$ are updated in each iteration based on the element-wise a-posteriori error indicators ${\eta_K}$. To support the adaptively distributed points selection, we introduce a sampling approach called recovery-type estimator based adaptive distribution (RecAD), which helps allocate adaptively distributed points to regions with larger estimated errors.

\begin{figure}[htbp]
    \centering
    \includegraphics[height=0.25\textheight]{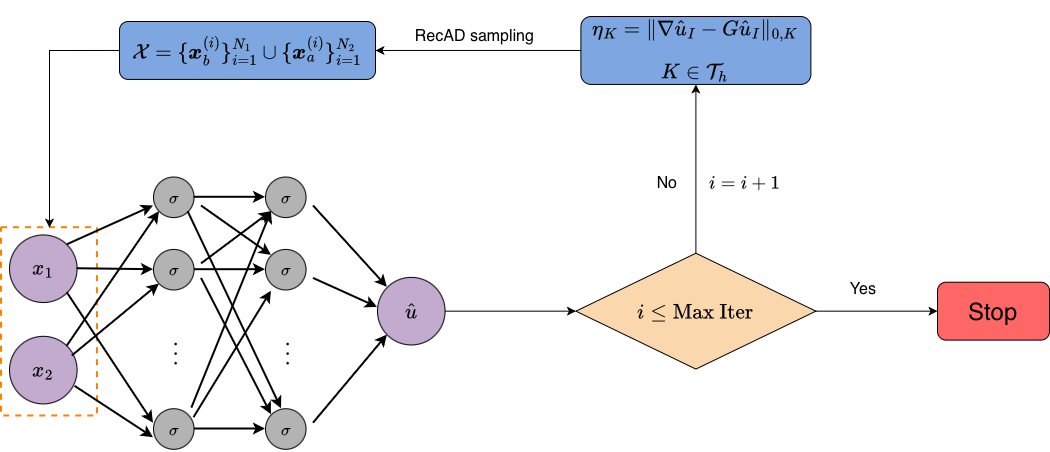}
    \caption{Flowchart of the R-PINN algorithm.}
    \label{flowchart}
\end{figure}

Unlike common sampling methods, the RecAD technique is element-based, as the computed a-posteriori errors are distributed across elements. This means we only need to specify the number of adaptively distributed points to sample in each element, allowing for random sampling within each element. Therefore, a reasonable scheme is required to assign the number of adaptively distributed points across elements based on their estimated a-posteriori errors. To this end, we design an iterative procedure.

We begin by sorting all elements in descending order according to their estimated a-posteriori errors, and denote them by the set $\mathcal{S}_0 = \{K_j\}_{j=1}^{N_e}$, where $N_e$ is the total number of elements. Let $n_1, n_2, \dots, n_{N_e}$ denote the number of adaptively distributed points assigned to each element, which are initially all set to zero.
Additionally, we introduce a tolerance parameter $\epsilon \in \left(\frac{1}{N_e}, 1\right)$, representing the proportion of elements to be selected for sampling. 

\begin{algorithm}[htbp]
	\caption{Recovery-type estimator based adaptive distribution (RecAD)}\label{read}
	\begin{algorithmic}[1] 
		
			\STATE Initialize the number of remaining adaptively distributed points $N_r = N_2$, the number of remaining elements $m_0 = N_e$, and set the tolerance parameter $\epsilon$.
			\STATE Sort all elements in descending order according to their estimated a-posteriori errors to obtain the set $\mathcal{S}_0 = \{K_j\}_{j=1}^{N_e}$.
			\STATE Set the initial number of assigned points for each element as $n_j = 0$, for $j=1,2,\cdots,N_e$.
			\STATE Set iteration counter $k \gets 1$.
		\WHILE{$N_r \ne 0$}
			\STATE Select $m_k = \left[ \epsilon m_{k-1} \right]$ elements with the largest a-posteriori errors from $\mathcal{S}_{k-1}$ as $\mathcal{S}_k$.
			\STATE Assign the number of adaptively distributed points to elements in $\mathcal{S}_k$ proportionally:
				$$n_j \gets n_j + \left[ N_r \frac{\eta_{K_j}}{\sum_{i=1}^{m_k} \eta_{K_i}} \right], \quad j = 1, 2, \cdots, m_k.$$
			\STATE Update the number of remaining adaptively distributed points $N_r$ by equation~(\ref{redundant}).
			\IF{$m_k<\frac{1}{\epsilon}$}
				\STATE Assign the remaining points to the element with the largest error: $n_1 \gets n_1 + N_r$.
				\STATE \textbf{break}
				\ENDIF
			\STATE $k\gets k+1$.
			\ENDWHILE
		\STATE Randomly sample the adaptively distributed points $\{\bm{x}_a\}_{i=1}^{N_2}$ with the numbers 
			$n_1, n_2,\cdots,n_{N_e}$ in each element.
		\STATE \textbf{Output:} Adaptively distributed points $\{\bm{x}_a\}_{i=1}^{N_2}$.
	\end{algorithmic}
\end{algorithm}

At the first iteration, we select the top $m_1 = \left[ \epsilon N_e \right]$ elements with the largest a-posteriori errors from $\mathcal{S}_0$, and denote this subset as $\mathcal{S}_1 = \{K_j\}_{j=1}^{m_1}$. In these selected elements, the number of adaptively distributed points are assigned approximately in proportion to their estimated a-posteriori errors,
\begin{equation*}
    n_j = \left[N_2\frac{\eta_{K_j}}{\sum_{i=1}^{m_1}\eta_{K_i}}\right],\quad j=1,2,\cdots,m_1.
\end{equation*}
Since rounding is applied to ensure that each $n_j$ is an integer, the total number of assigned points is likely to be smaller than the prescribed $N_2$. We then define the number of remaining points $N_r$ as
\begin{equation}\label{redundant}
    N_r=N_2-\sum_{j=1}^{N_e}n_j.
\end{equation}
In the next iteration, we select a subset of elements from $\mathcal{S}_1$, specifically the top $m_2 = \left[ \epsilon m_1 \right]$ elements with the largest a-posteriori errors, denoted as $\mathcal{S}_2 = \{K_j\}_{j=1}^{m_2}$. The number of remaining points $N_r$ are then redistributed to these elements according to their relative a-posteriori errors
\begin{equation*}
    n_j=n_j+\left[N_r\frac{\eta_{K_j}}{\sum_{i=1}^{m_2}\eta_{K_i}}\right],\quad j=1,2,\cdots,m_2.
\end{equation*}
After each iteration, we update $N_r$ using the same equation~(\ref{redundant}). The process continues until one of the following stopping criteria is met:
(1) the number of adaptively distributed points $N_2$ has been fully allocated (namely $N_r = 0$), or
(2) the number of selected elements becomes smaller than $\frac{1}{\epsilon}$, meaning the next subset would be empty. In the second case, suppose that after $k$-th iterations we obtain $\mathcal{S}_k = \{K_j\}_{j=1}^{m_k}$ with $N_r > 0$ and $m_k < \frac{1}{\epsilon}$. Then, continuing the iteration would result in an empty set $\mathcal{S}_{k+1}$. To avoid this, we directly assign the number of remaining points $N_r$ to the element with the largest a-posteriori error $K_1$, namely $n_1=n_1+N_r$.Now that the number of adaptively distributed points for each element in the mesh $\mathcal{T}_h$ has been fully determined, we then randomly sample the corresponding number of points within each element based on the assigned values $n_1, n_2, \dots, n_{N_e}$, thereby completing the adaptive sampling process.

The full algorithm for this adaptive distribution scheme RecAD is summarized in Algorithm~\ref{read}. With RecAD in place, all components of the R-PINN framework are complete. The entire training procedure of R-PINN is outlined in Algorithm~\ref{R-PINN}.

\begin{algorithm}[htbp]
	\caption{Recovery-type estimator enhanced PINN (R-PINN)}\label{R-PINN}
	\begin{algorithmic}[1] 
		\REQUIRE{Number of adaptively distributed points $N_a$, number of background points $N_b$, and maximum number of adaptive iterations $M$.}
		\STATE Generate $N_1$ background collocation points $\{\bm{x}_b^{(i)}\}_{i=1}^{N_1}$ with Sobol sequence. 
		\STATE Pre-train the network $\hat{u}_{\bm{\Theta}}(\bm{x})=\mathcal{N}^L(\bm{x};\bm{\Theta})$ with background collocation points $\{\bm{x}_b^{(i)}\}_{i=1}^{N_1}$.
		\STATE  Set the iteration counter $k \gets 1$.
            \WHILE{$k\leq M$}
            \STATE Compute the a-posteriori errors of $\hat{u}_{\bm{\Theta}}(\bm{x})$ with the recovery-type estimator defined in the equation~(\ref{estimator-eq}). 
            \STATE Generate $N_2$ adaptively distributed points $\{\bm{x}_a^{(i)}\}_{i=1}^{N_2}$ with RecAD method.
            \STATE Retrain the network $\hat{u}_{\bm{\Theta}}(\bm{x})=\mathcal{N}^L(\bm{x};\bm{\Theta})$ with the dataset $\mathcal{D}=\{\bm{x}_b^{(i)}\}_{i=1}^{N_1}\cup\{\bm{x}_a^{(i)}\}_{i=1}^{N_2}$.
            \ENDWHILE
		\RETURN  The approximated solution $\hat{u}_{\bm{\Theta}}(\bm{x})=\mathcal{N}^L(\bm{x};\bm{\Theta})$.
	\end{algorithmic}
\end{algorithm}

\section{Numerical experiments}

In this section, we present numerical experiments demonstrating the effectiveness of our proposed algorithm. Our evaluation compares three  kinds of PINNs: recovery-type estimator enhanced PINN (R-PINN), failure-informed PINN (FI-PINN) \cite{gao2023failure} and residual-based attention in PINN (RBA-PINN) \cite{ANAGNOSTOPOULOS2024116805}.  Among these methods, the first two adopt the adaptive sampling strategy, while RBA-PINN employs a strategy of adaptive weighting. For the gradient reconstruction in R-PINN, we employ the weighted averaging approach throughout all experiments. It is observed that the other two reconstruction methods achieve similarly comparable performance.

\subsection{Experimental setup}
To ensure a fair comparison, we maintain identical parameters and training processes across all methods for the same number of iterations, varying only the number and distribution of adaptively distributed points. The recovery-type a-posteriori estimation is implemented using \textit{FEALPy} \cite{fealpy}.

All experiments employ a fully-connected neural network with $7$ hidden layers, each containing 20 neurons with \textit{tanh} for the activation function unless otherwise specified. The training process comprises two phases: an initial pre-training phase utilizing the \textit{L-BFGS} optimizer with a learning rate of $0.1$, conducted for a problem-dependent number of epochs, followed by an adaptive training phase consisting of $M$ adaptive iterations. Each iteration is constrained to $5000$ epochs while employing the same optimization method unless otherwise specified.
\begin{figure}
    \centering
    \includegraphics[height=0.2\textheight]{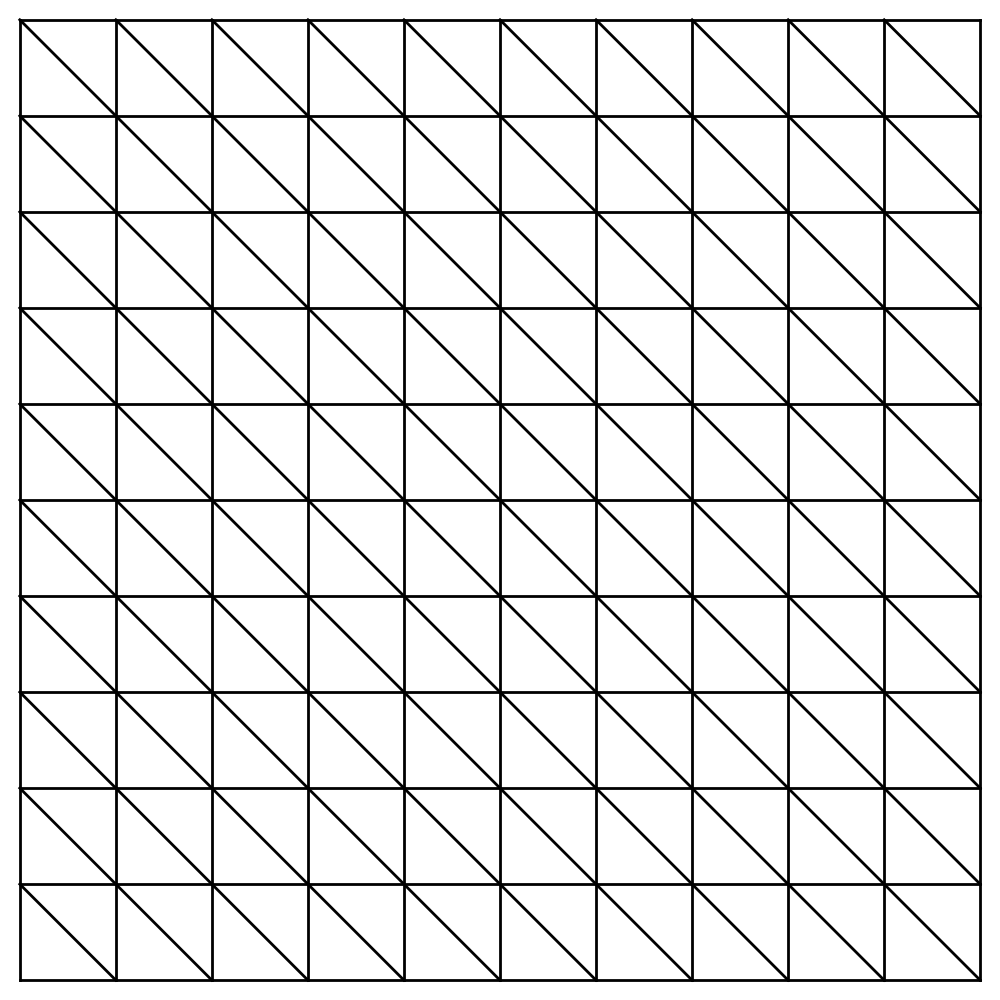}
    \caption{Illustration of a $10\times10$ regular diagonal triangulation.}
    \label{mesh}
\end{figure}

For the choice of parameters, we set $\epsilon = 0.02$ for R-PINN. The triangulation of the domain adopted in the method is a uniform diagonal mesh composed of right triangles. An example with a $10\times 10$ mesh is illustrated in Fig.~\ref{mesh}. Unless otherwise specified, we use a $50\times 50$ mesh as the default configuration in the following experiments.
 For collocation points, we use $N_1 = 2000$ background collocation points and $N_2 = 100$ adaptively distributed points for R-PINN unless specified. For the RBA-PINN, only $N_1=2000$ background collocation points are used during training, without incorporating additional adaptively distributed points unless specified. However, the training process remains the same as that of the two adaptive sampling PINNs. All methods utilize $200$ boundary collocation points. Notably, our approach achieves comparable or superior accuracy despite using fewer total epochs and adaptively distributed points than related studies. The larger $N_2$ value for FI-PINN is necessary as it fails to achieve adaptivity with fewer points. Furthermore, while our method maintains a fixed number of adaptively distributed points, FI-PINN accumulates points throughout the adaptive iteration process.
To evaluate the accuracy of our method, we employ two error metrics: the relative $L^2$ error,

\begin{equation*}
	\frac{\lVert\hat{u}_{\bm{\Theta}}-u\rVert_2}{\lVert u\rVert_2},
\end{equation*}
and the $L^{\infty}$ error, defined as $\lVert\hat{u}_{\bm{\Theta}} - u\rVert_\infty$. Here, $\hat{u}_{\bm{\Theta}}$ denotes the numerical solution obtained by our method, and $u$ represents the exact solution.

\subsection{Poisson's equation\label{poisson}}
Consider the following Poisson's equation:
\begin{equation*}
	\begin{split}
		-\Delta u(x,y)=&f(x,y),\quad (x,y)\ \text{in}\ \Omega, \\
		u(x,y) = &g(x,y), \quad (x,y)\ \text{on}\ \partial\Omega,
	\end{split}
\end{equation*}
where the domain $\Omega=[0,1]\times[0,1]$, and the true solution is 
\begin{equation*}
	u(x,y)=e^{-1000[(x-0.5)^2+(y-0.5)^2]}.
\end{equation*}
The solution exhibits a sharp peak near the point $(0.5, 0.5)$, which is challenging for vanilla PINN to capture accurately. In this experiment, the network is pre-trained for 15,000 epochs. The number of background collocation points is set to $N_1=4000$, and the maximum number of adaptive iterations is specified as $M=4$. Each adaptive iteration consists of 10,000 epochs instead of 5,000, due to the difficulty of approximating this particular example, which features a solution with a sharp peak.

\begin{table}[htbp]
    \centering
    \begin{tabular}{||c||ccc||ccc||cc||}
        \hline
       & \multicolumn{3}{c||}{R-PINN} & \multicolumn{3}{c||}{FI-PINN} & \multicolumn{2}{c||}{RBA-PINN} \\
        \hline
        $N_{\text{iter}}$ & Relative $L^2$ & $L^\infty$ &$N_2$  & Relative $L^2$ & $L^\infty$ & $N_2$& Relative $L^2$ & $L^\infty$\\
        \hline\hline
        1    & 0.1695 & 0.0259 & 120 & 0.4166 & 0.2862 & 500 & 0.3042 & 0.2040 \\
        2    & 0.0673 & 0.0177 & 120 & 
        0.1656 & 0.0449 & 1000 & 0.3360  & 0.1880 \\
        3    & 0.0497 & 0.0103 & 120 & 
        0.1405 & 0.0160 & 1500 & 0.2781 & 0.1828 \\
        4    & 0.0277 & 0.0015  & 120&
        0.0530 & 0.0086 & 2000 & 0.2741 & 0.1791 \\
        \hline
    \end{tabular}
    \caption{The relative $L^2$ and $L^\infty$ error obtained by three different PINNs solving Poisson's equation.}
    \label{Poisson-err}
\end{table}

\begin{figure}[t]
	\centering
	\subcaptionbox*{R-PINN}
	{
		\includegraphics[height=0.2\textheight]{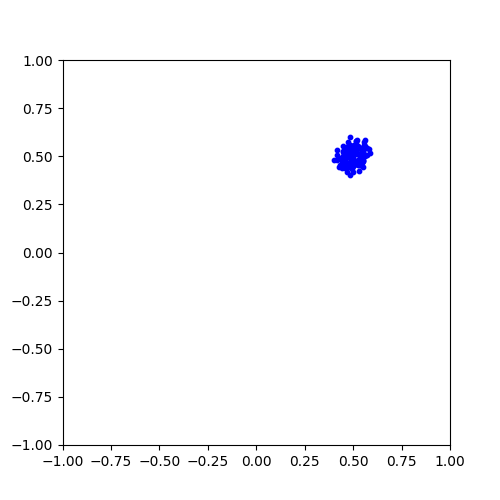}
		\includegraphics[height=0.2\textheight]{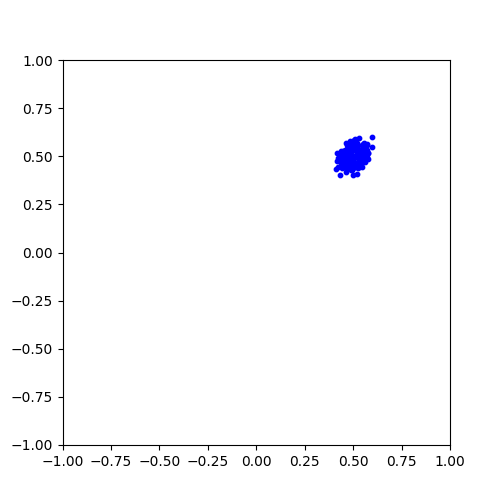}
	}
	\subcaptionbox*{FI-PINN}
	{
		\includegraphics[height=0.2\textheight]{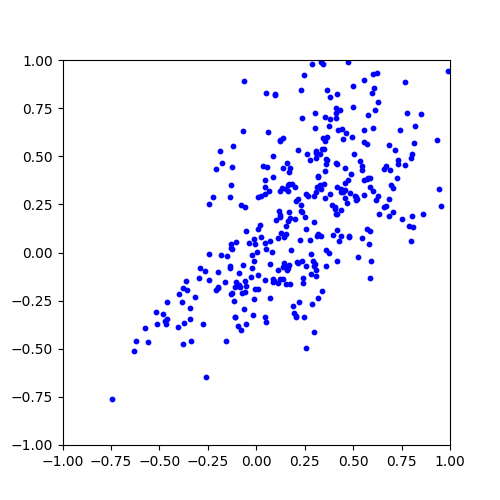}
		\includegraphics[height=0.2\textheight]{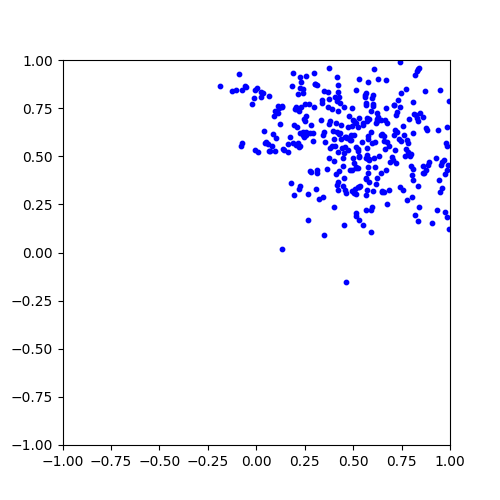}
	}
	\caption{The distribution of adaptively distributed points in the first $2$ adaptive iterations with two adaptive PINN methods.}
	\label{Poisson-points}
\end{figure}
\begin{figure}[p]
	\centering
	\subcaptionbox*{R-PINN}
	{
        \includegraphics[height=0.2\textheight]{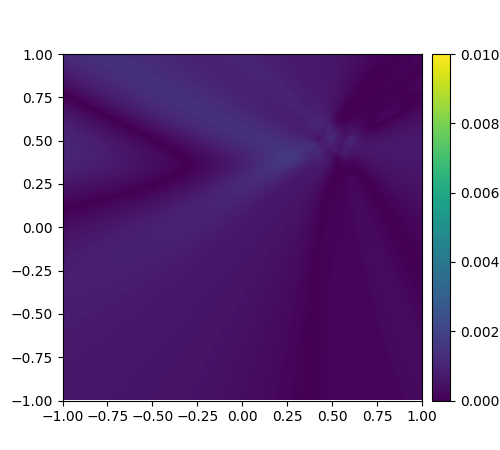}
		\includegraphics[height=0.2\textheight]{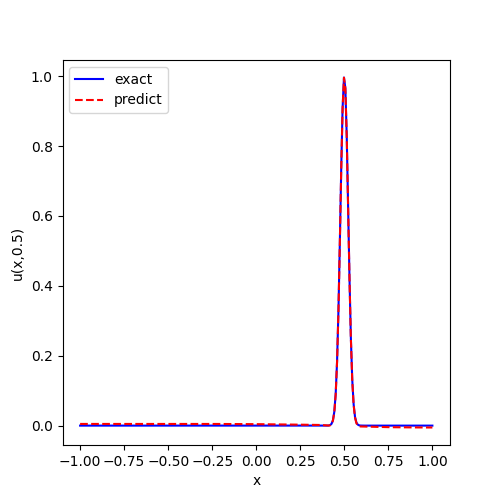}
		\includegraphics[height=0.2\textheight]{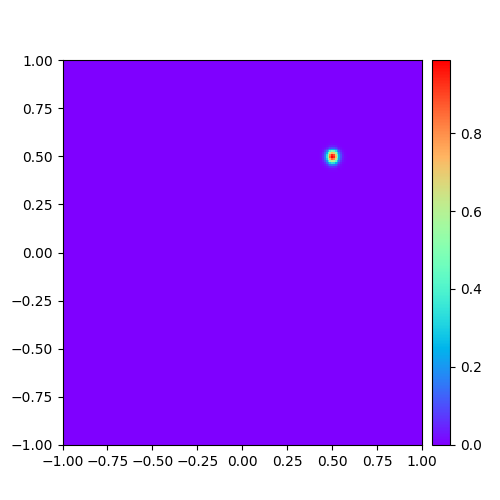}
	}
        \subcaptionbox*{FI-PINN}
	{
        \includegraphics[height=0.2\textheight]{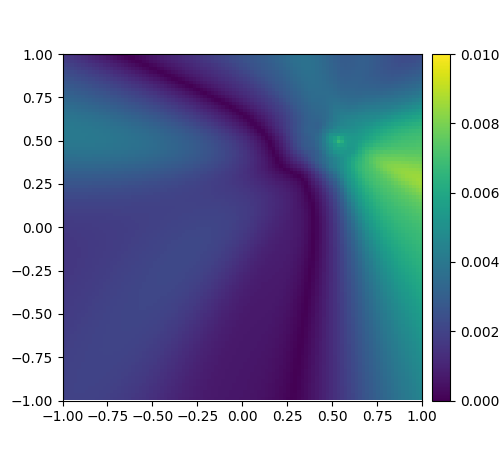}
		\includegraphics[height=0.2\textheight]{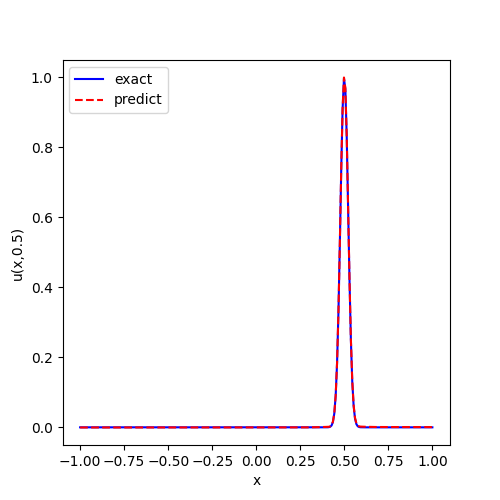}
		\includegraphics[height=0.2\textheight]{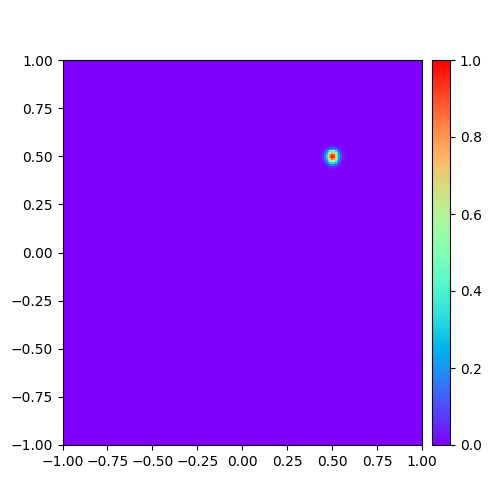}
	}
	\subcaptionbox*{RBA-PINN}
	{
        \includegraphics[height=0.2\textheight]{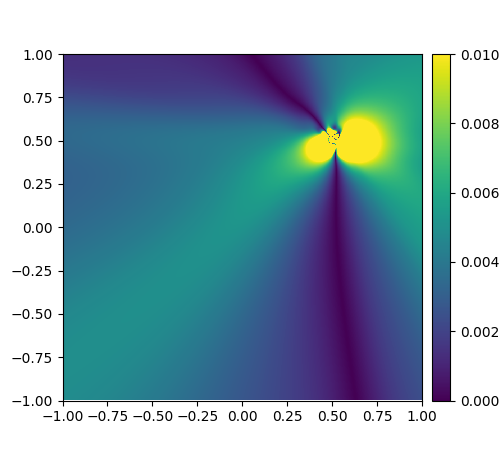}
		\includegraphics[height=0.2\textheight]{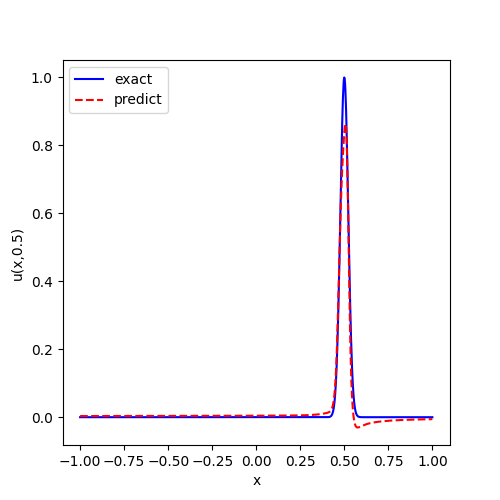}
		\includegraphics[height=0.2\textheight]{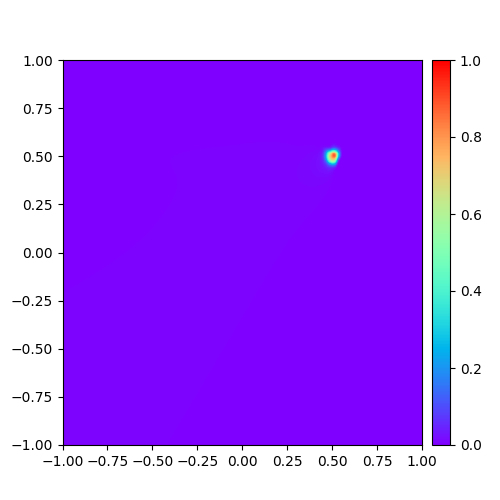}
	}
	\caption{The absolute errors across the domain (left), the cross-sections at $y=0.5$ (middle) and the predicted solutions of three adaptive methods solving Poisson's equation.}
	\label{Poisson-err-curve}
\end{figure}
\begin{figure}[htbp]
    \centering
    \includegraphics[height=0.28\textheight]{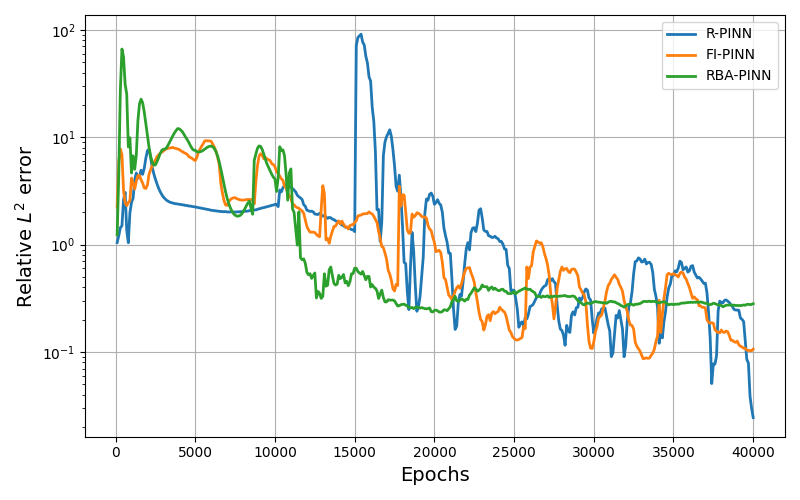}
    \caption{The relative $L^2$ error by three different PINNs in the optimization process of solving Poisson's equation.}
    \label{Poisson-L2}
\end{figure}

Fig.~\ref{Poisson-points} illustrates the distribution of adaptively distributed points by two adaptive sampling methods in the first two adaptive iterations. It is evident that all adaptively distributed points of our model converge to the neighborhood of the point $(0.5,0.5)$, while adaptively distributed points of FI-PINN are more dispersed. This shows the strong effectiveness of our algorithm in handling problems with significant local variations. Table~\ref{Poisson-err} compares the relative $L^2$ error and $L^\infty$ error among R-PINN, FI-PINN and RBA-PINN during the adaptive iterations where $N_{\text{iter}}$ represents the number of adaptive iterations. The relative $L^2$ error for R-PINN reaches 0.0277, and the $L^\infty$ error reaches 0.0015.
In contrast, the errors obtained by FI-PINN decrease slowly at the beginning, but drop fast later, reaching the relative $L^2$ error of 0.0530 and $L^\infty$ error of 0.0086. It can be observed that R-PINN can converge faster than FI-PINN. R-PINN can also  reach a smaller relative $L^2$ error even in the case where the adaptively distributed points of FI-PINN are continuously accumulated. As for RBA-PINN, it appears to have a difficulty in capturing the local features of the solution accurately. Similar results are observed in Fig.~\ref{Poisson-L2}, where the lowest relative $L^2$ error for our algorithm is below $5\times 10^{-2}$. In the left panel of Fig.~\ref{Poisson-err-curve}, we show the absolute errors throughout the entire domain, clearly demonstrating that the errors are much smaller with R-PINN. The middle panel of Fig.~\ref{Poisson-err-curve} illustrates the cross-sections at $y=0.5$. Both R-PINN and FI-PINN are able to accurately approximate the peak-shaped solution. In contrast, although RBA-PINN identifies the local feature, its accuracy is limited by the insufficient number of collocation points in the corresponding region. The final solutions by three different PINNs are plotted in the right of Fig.~\ref{Poisson-err-curve}.

\subsection{Burgers' equation}

We consider the following Burgers' equation:
\begin{equation*}
	\begin{split}
		&\frac{\partial u}{\partial t}+u\frac{\partial u}{\partial x}=\frac{0.001}{\pi}\frac{\partial^2 u}{\partial x^2},\\
		&u(x,0)=-\sin(\pi x),\\
		&u(t,-1)=u(t,1)=0
	\end{split}
\end{equation*}
where the domain is $(x,t)\in [-1,1]\times [0,1]$. The solution to this Burgers' equation exhibits an extremely sharp mutation along $x=0$, which makes the problem highly susceptible to numerical errors and local deviations. In this test, we pre-train the network with $N_1=50,000$ background collocation points for $15,000$ iterations and use $M=4$ steps of adaptive processes.

\begin{table}[htbp]
    \centering
    \begin{tabular}{||c||ccc||ccc||cc||}
        \hline
        & \multicolumn{3}{c||}{R-PINN} & \multicolumn{3}{c||}{FI-PINN} & \multicolumn{2}{c||}{RBA-PINN} \\
        \hline
        $N_{\text{iter}}$ & Relative $L^2$ & $L^\infty$ & $N_2$ & Relative $L^2$ & $L^\infty$ & $N_2$ & Relative $L^2$ & $L^\infty$  \\
        \hline\hline
        1    & 0.0166 & 0.3260 & 500 & 0.2106 & 1.8845 & 600 & 0.1692 & 1.9880 \\
        2    & 0.0169 & 0.4087 & 500 & 0.1997 & 1.9430 & 1200 & 0.1692 & 1.9881 \\
        3    & 0.0122 & 0.2078 & 500 & 0.1926 & 1.9593 & 1800 & 0.1693 & 1.9881 \\
        4    & 0.0048 & 0.0920 & 500 & 0.1928 & 1.9471 & 2400 & 0.1693 & 1.9881 \\
        \hline
    \end{tabular}
    \caption{The relative $L^2$ and $L^\infty$ error obtained by three different PINNs solving Burgers' equation.}
    \label{Burgers-err}
\end{table}

\begin{figure}[htbp]
    \centering
    \subcaptionbox*{R-PINN}
    {
    \includegraphics[height=0.16\textheight]{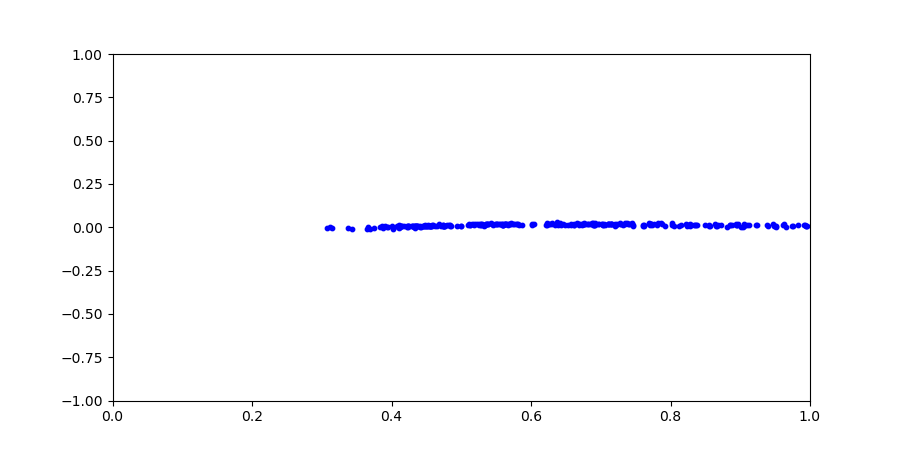}
    \includegraphics[height=0.16\textheight]{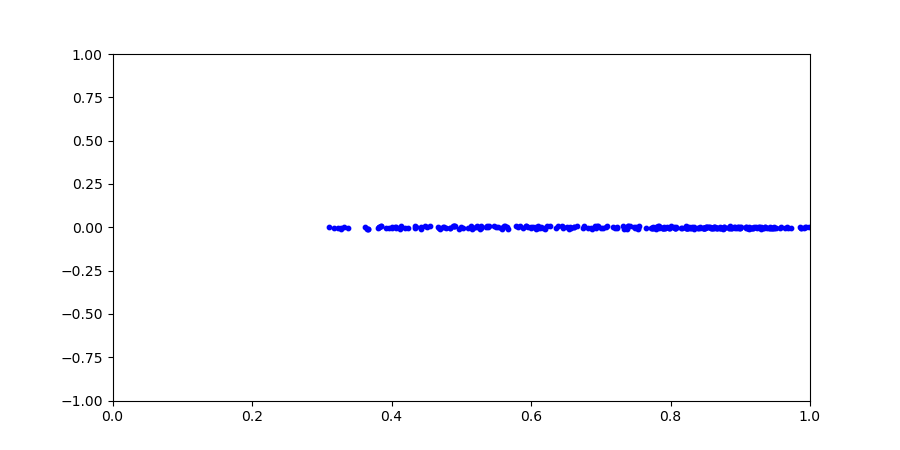}
	}
	\subcaptionbox*{FI-PINN}
	{
    \includegraphics[height=0.16\textheight]{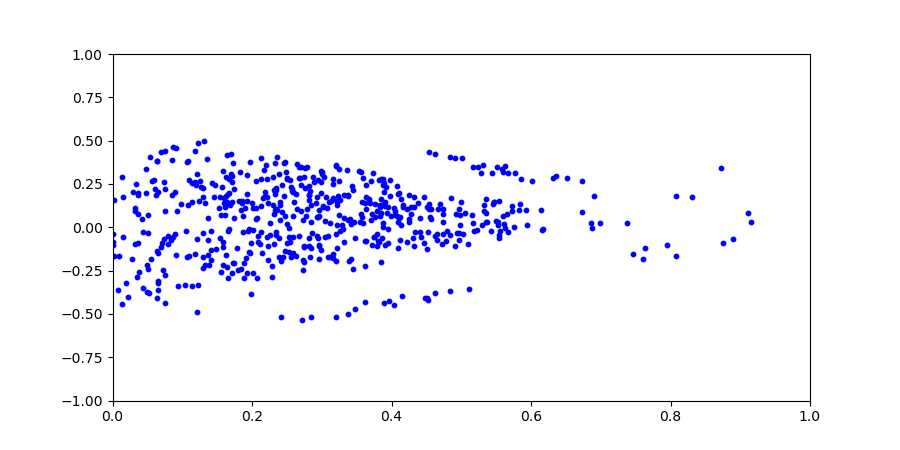}
    \includegraphics[height=0.16\textheight]{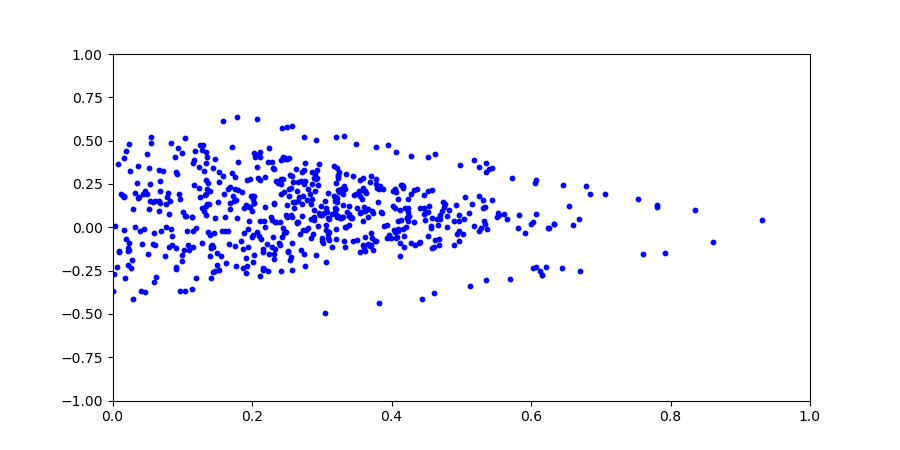}
	}
    \caption{The distribution of adaptively distributed points in the first $2$ adaptive iterations with two adaptive PINN methods solving Burgers' equation.}
    \label{Burgers-points}
\end{figure}

\begin{figure}[p]
	\centering
	\subcaptionbox*{R-PINN}
	{
        \includegraphics[height=0.16\textheight]{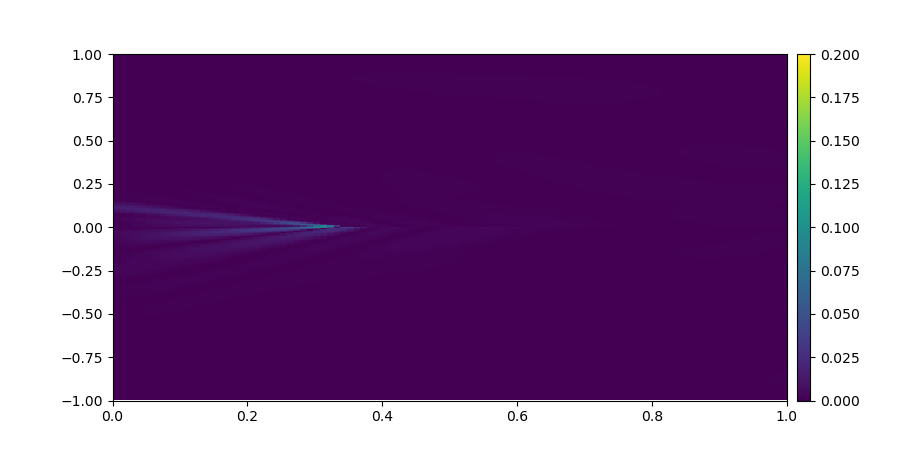}
		\includegraphics[height=0.16\textheight]{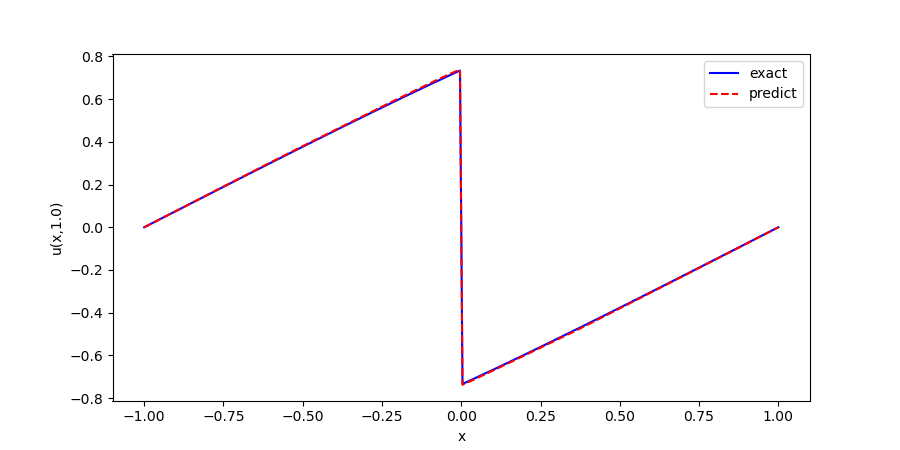}
	}
	\subcaptionbox*{FI-PINN}
	{
		\includegraphics[height=0.16\textheight]{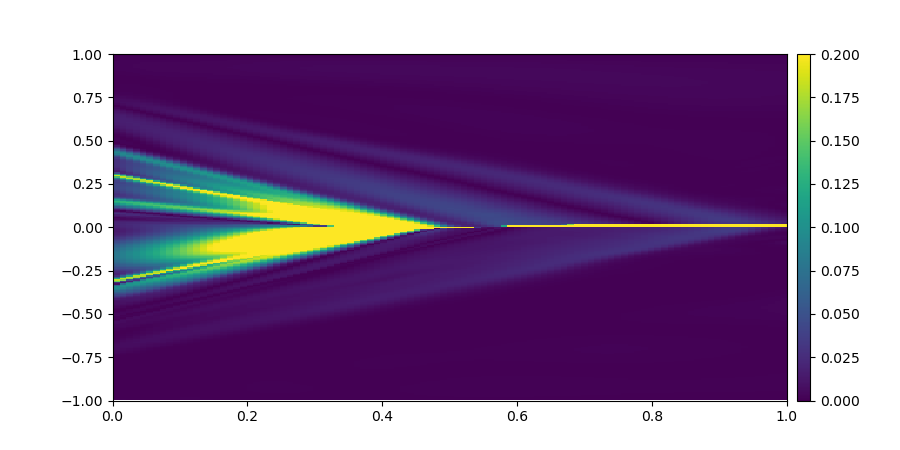}
		\includegraphics[height=0.16\textheight]{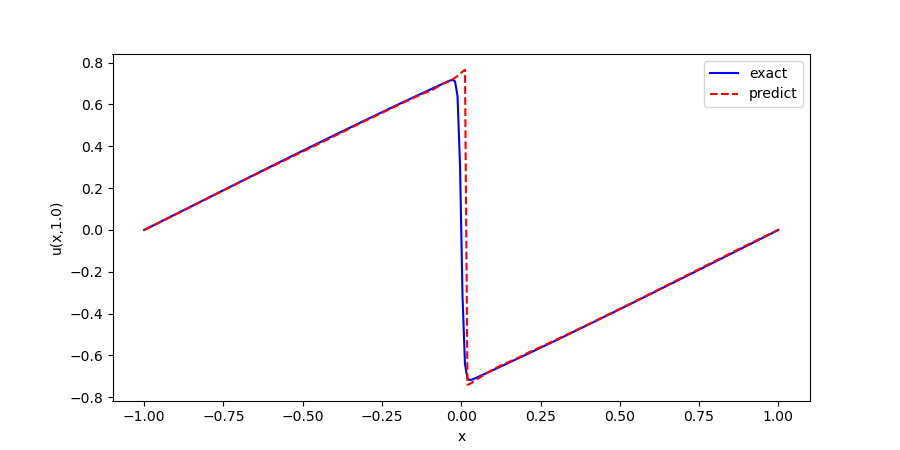}
	}
	\subcaptionbox*{RBA-PINN}
	{
    \includegraphics[height=0.16\textheight]{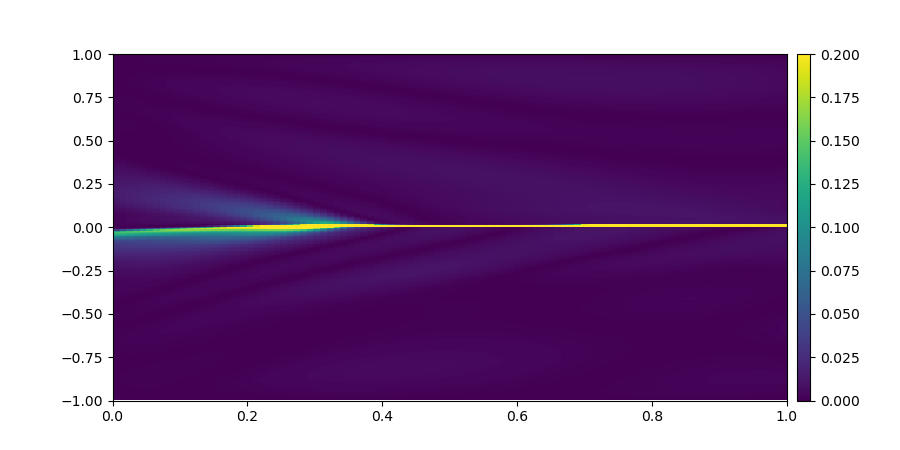}
		\includegraphics[height=0.16\textheight]{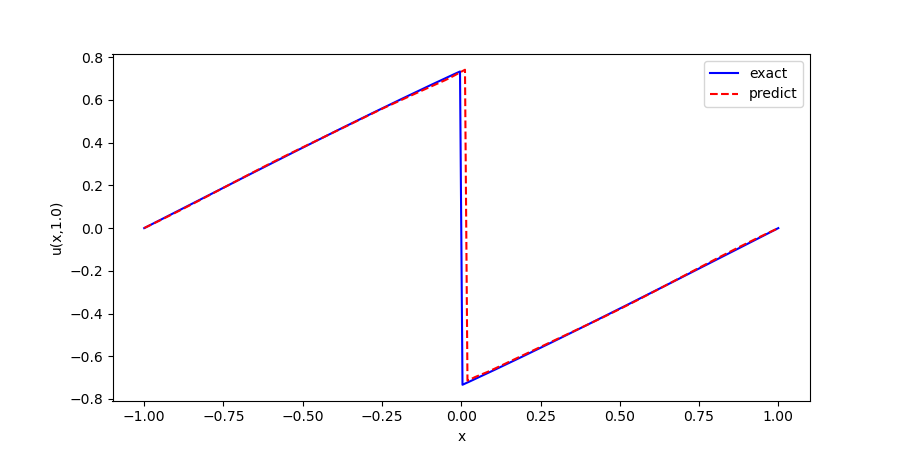}
	}
	\caption{The absolute errors across the domain (left) and the cross-sections at $t=1$ (right) of three different PINNs solving Burgers' equation.}
	\label{Burgers-err-curve}
\end{figure}

\begin{figure}[htbp]
    \centering
    \includegraphics[height=0.16\textheight]{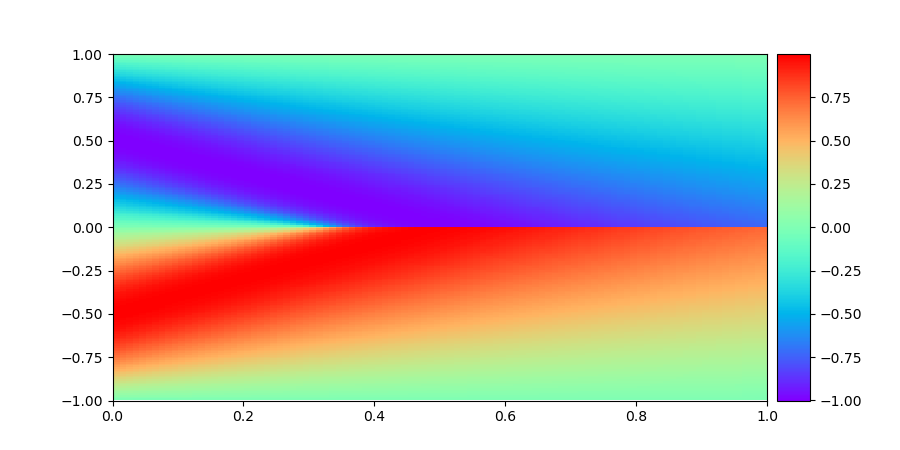}
    \includegraphics[height=0.16\textheight]{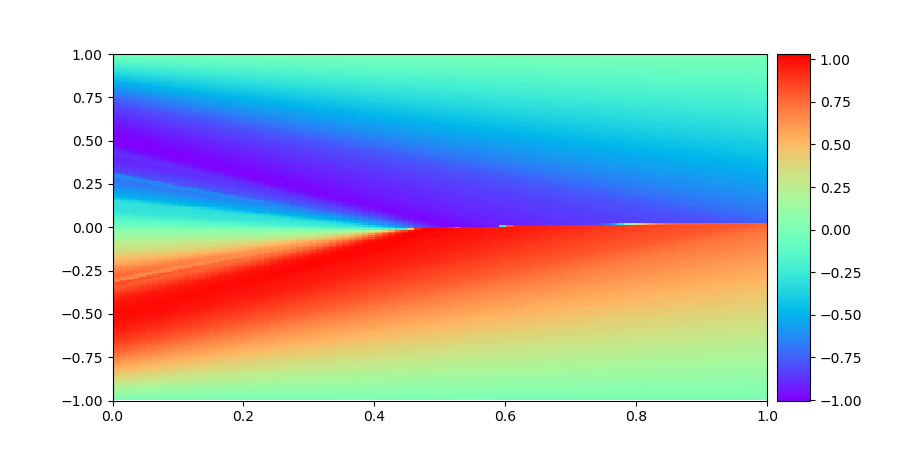}
    \caption{The predicted solution (R-PINN at the left and FI-PINN at the right).}
    \label{Burgers-sol}
\end{figure}

\begin{figure}[htbp]
    \centering
    \includegraphics[height=0.28\textheight]{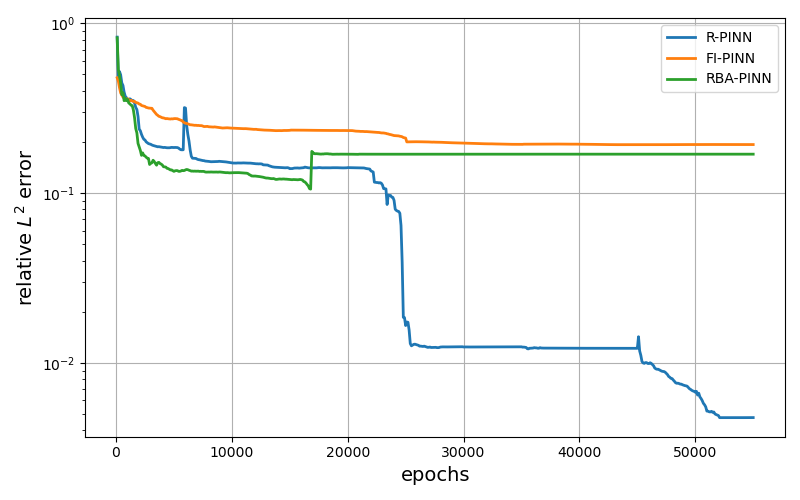}
    \caption{The relative $L^2$ error by three different PINNs in the optimization process of solving Burgers' equation.}
    \label{Burgers-L2}
\end{figure}

The adaptively distributed points generated during the first two adaptive training processes by two adaptive sampling PINN methods are illustrated in Fig.~\ref{Burgers-points}. We observe that the adaptively distributed points successfully cluster around $x=0$ with both methods, where the error is largest (as seen in Fig.~\ref{Burgers-err-curve}). This clustering significantly improves the local error. The exact relative $L^2$ error and $L^\infty$ error are listed in Table~\ref{Burgers-err}. As shown, the final relative $L^2$ error and $L^\infty$ error for R-PINN are 0.0048 and 0.0920, respectively, while those for FI-PINN are 0.1928 and 1.9471. In contrast, the final relative $L^2$ and $L^\infty$ errors of RBA-PINN are 0.1693 and 1.9881, respectively. Both FI-PINN and RBA-PINN fail to accurately capture the shock, whereas R-PINN exhibits significantly lower errors than the other two methods. Such behaviour can also be observed in Fig.~\ref{Burgers-L2}, where the relative $L^2$ error is illustrated specifically, showing that the lowest error can reach below $1\times 10^{-2}$. We also show the absolute errors across the entire domain and the cross-sections at $t=1$ in Fig.~\ref{Burgers-err-curve}. The predicted solutions of two adaptive sampling methods are illustrated in Fig.~\ref{Burgers-sol}.

\subsection{Two-dimensional problem with two peaks}
Considering the following two-dimensional problem
\begin{equation*}
    \begin{split}
        -\nabla\cdot [u(x,y)\nabla (x^2+y^2)]+\nabla^2u(x,y)=f(x,y),\quad (x,y)\in\Omega,\\
        u(x,y)=b(x,y),\quad (x,y)\in\partial\Omega,
    \end{split}
\end{equation*}
where $\Omega=[-1,1]\times[-1,1]$. The true solution is
\begin{equation*}
    u(x,y)=e^{-1000[(x-0.5)^2+(y-0.5)^2]}+e^{-1000[(x+0.5)^2+(y+0.5)^2]},
\end{equation*}
which has two peaks at $(-0.5, -0.5)$ and $(0.5,0.5)$. This test mainly can show the ability of R-PINN to deal with multiple peaks. For this experiment, we pre-train our network for the period of $10000$ iterations and set the maximum number of adaptive iterations as $M=4$. Since two peaks problem is challenging for neural network to approximate, we use $5$ hidden layers with $64$ neurons in each layer to enhance its capability. Moreover, we use $N_2=150$ adaptively distributed points for R-PINN, and $N_2=750$ added adaptively distributed points for FI-PINN.

\begin{table}[htbp]
    \centering
    \begin{tabular}{||c||ccc||ccc||cc||}
        \hline
        & \multicolumn{3}{c||}{R-PINN} & \multicolumn{3}{c||}{FI-PINN} & \multicolumn{2}{c||}{RBA-PINN} \\
        \hline
        $N_{\text{iter}}$ & Relative $L^2$ & $L^\infty$ & $N_2$ & Relative $L^2$ &$L^\infty$ & $N_2$ & Relative $L^2$ & $L^\infty$ \\
        \hline\hline
        1    & 0.6778 & 0.1228 & 150 & 0.4372 & 0.4845 & 750 & 0.4045 & 0.5299 \\
        2    & 0.1202 & 0.0198 & 150 & 0.5220 & 0.3429 & 1500 & 0.3764 & 0.4885 \\
        3    & 0.1080 & 0.0259 & 150 & 0.4773 & 0.2233 & 2250 & 0.3578 & 0.4646 \\
        4    & 0.0840 & 0.0127 & 150 & 0.3969 & 0.2263 & 3000 & 0.3474 & 0.4507 \\
        \hline
    \end{tabular}
    \caption{The relative $L^2$ and $L^\infty$ error obtained by three different PINNs solving the two-dimensional problem with two peaks.}
    \label{twopeaks-err}
\end{table}

\begin{figure}[htbp]
    \centering
    \subcaptionbox*{R-PINN}
    {
    \includegraphics[height=0.2\textheight]{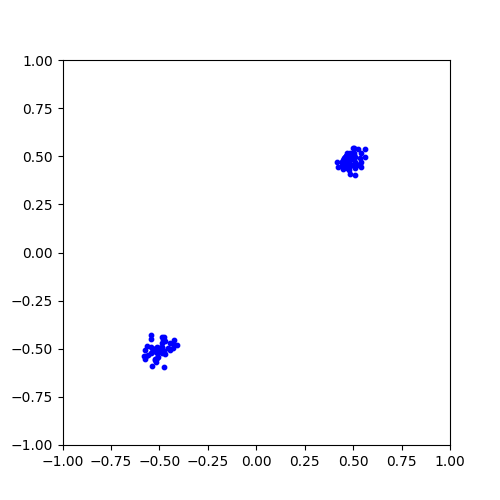}
    \includegraphics[height=0.2\textheight]{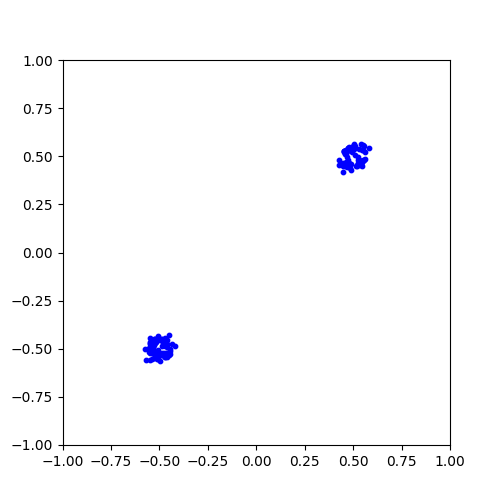}
	}
	\subcaptionbox*{FI-PINN}
	{
    \includegraphics[height=0.2\textheight]{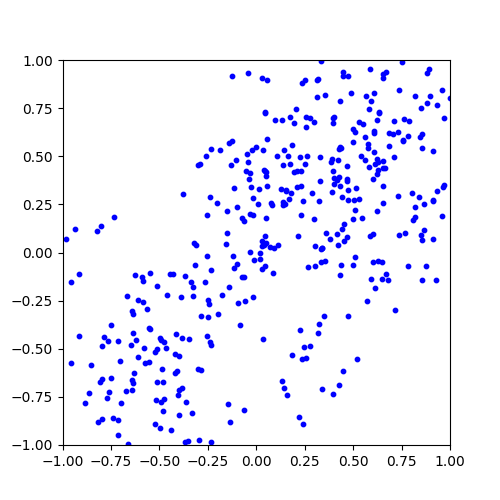}
    \includegraphics[height=0.2\textheight]{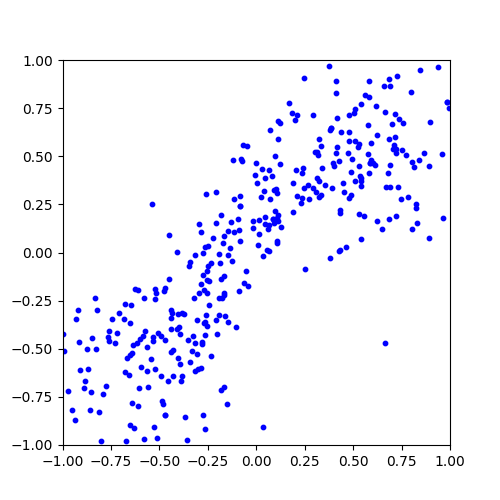}
	}
    \caption{The distribution of adaptively distributed points in the first two adaptive iterations of two adaptive sampling methods solving the two-dimensional problem with two peaks.}
    \label{twopeaks-points}
\end{figure}

\begin{figure}
	\centering
	\subcaptionbox*{R-PINN}
	{
        \includegraphics[height=0.2\textheight]{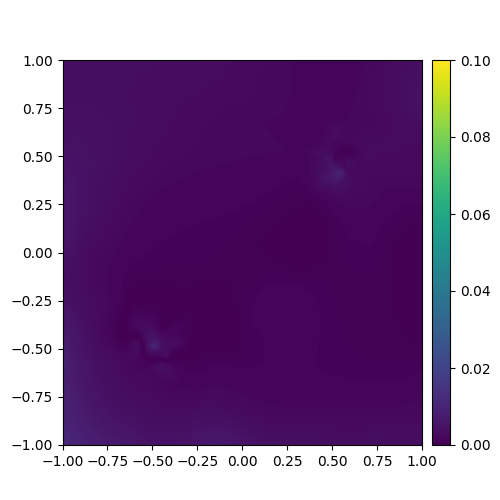}
        \includegraphics[height=0.2\textheight]{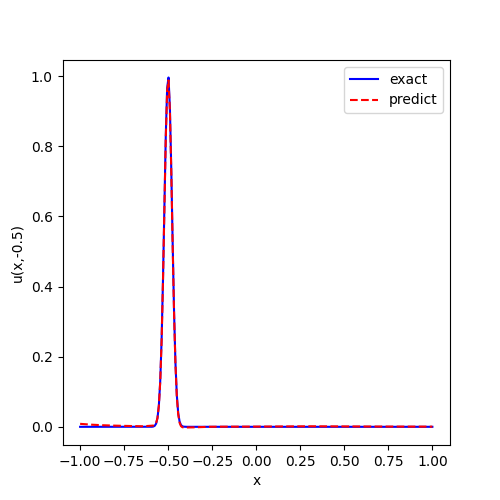}

        \includegraphics[height=0.2\textheight]{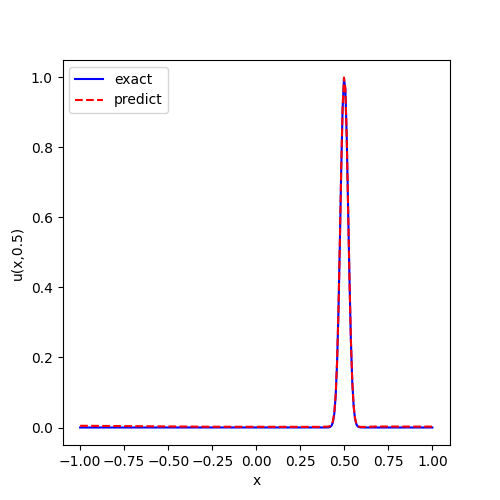}
	}
	\subcaptionbox*{FI-PINN}
	{
		\includegraphics[height=0.2\textheight]{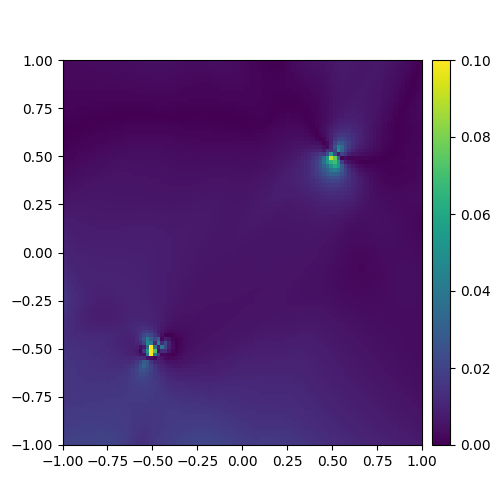}
		\includegraphics[height=0.2\textheight]{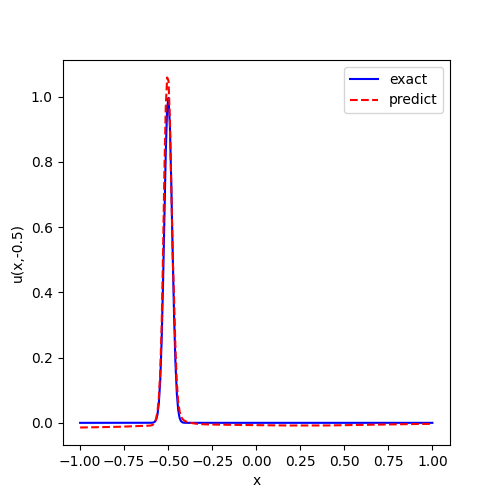}
		
		\includegraphics[height=0.2\textheight]{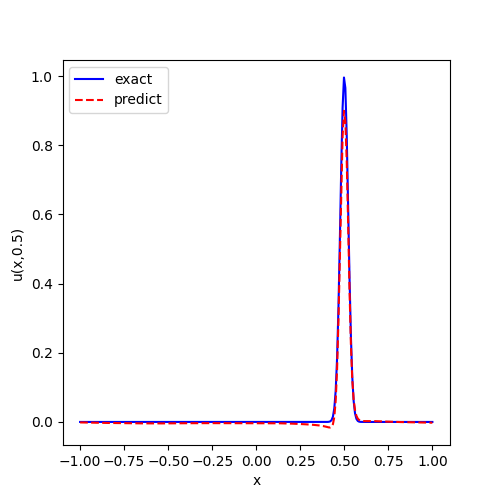}
	}
	\subcaptionbox*{RBA-PINN}
	{
        \includegraphics[height=0.2\textheight]{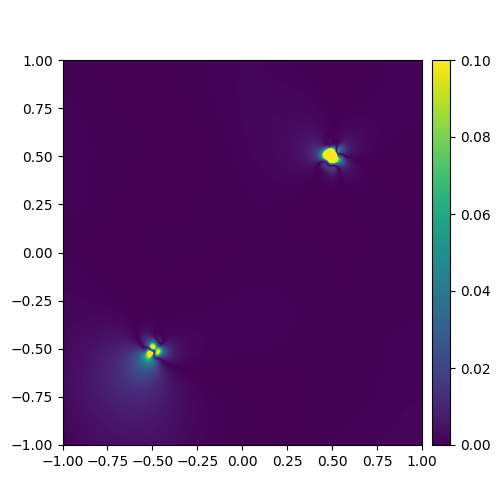}
		\includegraphics[height=0.2\textheight]{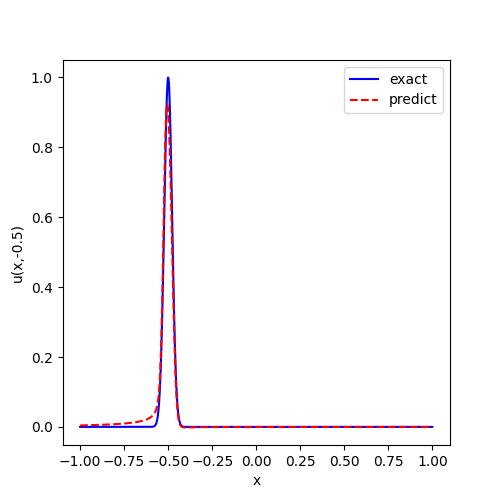}

        \includegraphics[height=0.2\textheight]{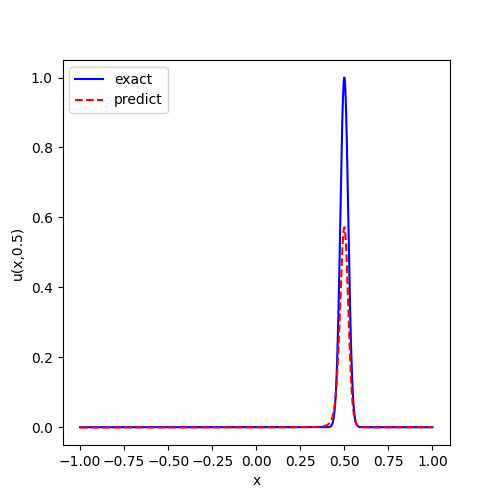}
	}
	\caption{The absolute errors across the domain (left), the cross-sections at $y=-0.5$ (middle) and the cross-sections at $y=0.5$ (right) of three different PINNs solving the two-dimensional problem with two peaks.}
	\label{twopeaks-err-curve}
\end{figure}

\begin{figure}[htbp]
    \centering
    \includegraphics[height=0.2\textheight]{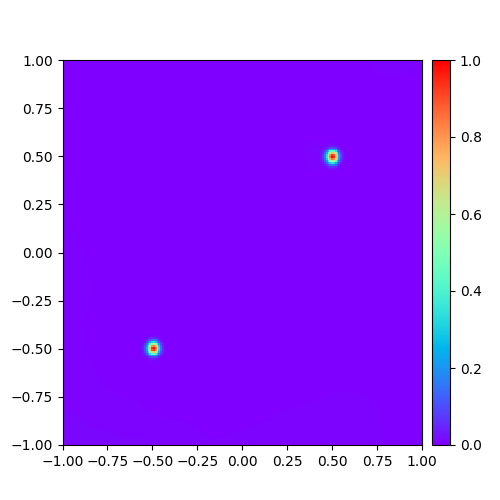}
    \includegraphics[height=0.2\textheight]{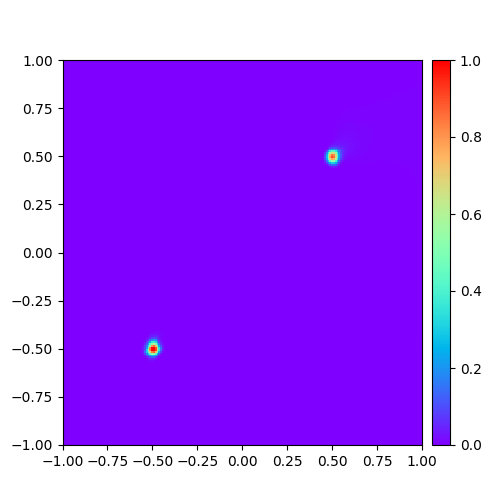}
    \includegraphics[height=0.2\textheight]{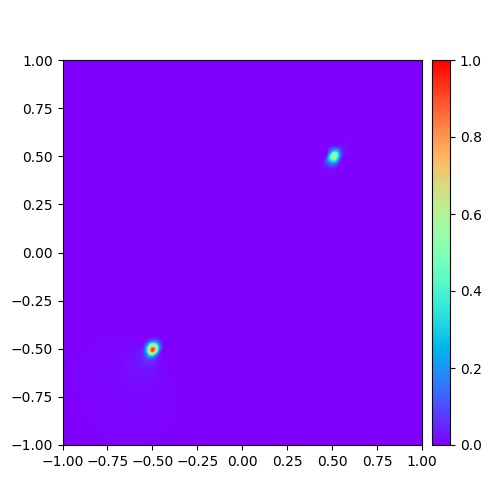}
    \caption{The predicted solution (R-PINN at the left, FI-PINN at the middle and RBA-PINN at the right) of the two-dimensional problem with two peaks.}
    \label{twopeaks-sol}
\end{figure}

\begin{figure}[htbp]
    \centering
    \includegraphics[height=0.28\textheight]{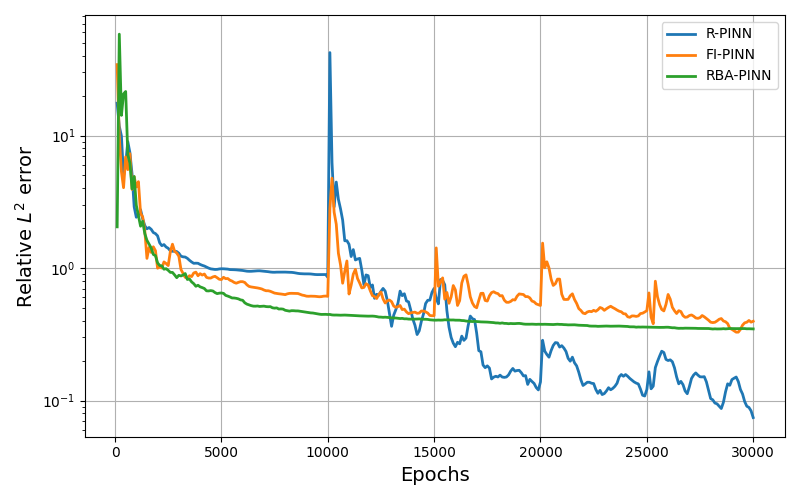}
    \caption{The relative $L^2$ error by three different PINNs in the optimization process of solving the two-dimensional problem with two peaks.}
    \label{twopeaks-L2}
\end{figure}

We present the distribution of adaptively distributed points of two adaptive sampling methods across the first two adaptive iterations in Fig.~\ref{twopeaks-points}. It is evident that these points converge around the vicinities of $(0.5, 0.5)$ and $(-0.5, -0.5)$ with R-PINN, demonstrating the capability of our algorithm to effectively handle problems characterized by multiple regions of abrupt changes. However, the distribution of adaptively distributed points with FI-PINN seems not so effective to deal with the two peaks. The reason of it might be a truncated Gaussian distribution is used, so that problems of more than one regions with large error are difficult to handle. This problem is also reflected in the training process. Table~\ref{twopeaks-err} compares the relative $L^2$ and $L^\infty$ errors among R-PINN, FI-PINN and RBA-PINN during the adaptive iterations. For R-PINN, the relative $L^2$ error reduces to 0.0840, and the $L^\infty$ error decreases to 0.0127. In contrast, the relative $L^2$ error associated with FI-PINN can only reach 0.3969 and the $L^\infty$ error can get to 0.2263, while RBA-PINN achieves a lower relative $L^2$ error 0.3474 but a larger $L^\infty$ error 0.4507. Similar trends are observed in Fig.~\ref{twopeaks-L2}, where the lowest relative $L^2$ error for R-PINN falls below $1 \times 10^{-1}$.

Fig.~\ref{twopeaks-err-curve} provides additional insights. The left panel illustrates the absolute errors across the entire domain, revealing significantly smaller errors for R-PINN. The middle and right panels display cross-sectional views at $y = -0.5$ and $y = 0.5$, respectively, showing that the predictions from R-PINN align closely with the exact solution, whereas there are large errors near two peaks for FI-PINN and RBA-PINN. Finally, the numerical solutions obtained by the three methods are visualized in Fig.~\ref{twopeaks-sol}.

\begin{table}[htbp]
	\centering
	\begin{tabular}{||c||ccc||ccc||cc||}
		\hline
		& \multicolumn{3}{c||}{R-PINN} & \multicolumn{3}{c||}{FI-PINN} & \multicolumn{2}{c||}{RBA-PINN} \\
		\hline
		$N_{\text{iter}}$ & Relative $L^2$ & $L^\infty$ & $N_2$ & Relative $L^2$ &$L^\infty$ & $N_2$ & Relative $L^2$ & $L^\infty$ \\
		\hline\hline
		1    & 0.0572 & 0.0761 & 100 & 0.1110 & 0.1590 & 500 & 0.1153 & 0.0458 \\
		2    & 0.0076 & 0.0098 & 100 & 0.0317 & 0.0491 & 1000 & 0.1140 & 0.0453 \\
        3 & 0.0076 & 0.0099 & 100 &
        0.0317 & 0.0490 & 1500 & 0.1137 & 0.0452\\
        4 & 0.0076 & 0.0099 & 100 & 0.0270 & 0.0418 & 2000 & 0.1135 & 0.0452\\
		\hline
	\end{tabular}
	\caption{The relative $L^2$ and $L^\infty$ error obtained by three different PINNs solving the wave equation.}
	\label{timedependent-l2}
\end{table}

\subsection{Wave equation}
For this test, we consider a wave equation as follows,
\begin{equation*}
    \begin{split}
        &\frac{\partial^2 u}{\partial t^2}-3\frac{\partial^2u}{\partial x^2}=0,\quad (t,x)\in [0,6]\times [-5,5],\\
        &u(0,x)=\frac{1}{\cosh(2x)}-\frac{0.5}{\cosh(2(x-10))}-\frac{0.5}{\cosh(2(x+10))},\\
        &\frac{\partial u}{\partial t}(0,x)=0,\\
        &u(t,-5)=u(t,5)=0,
    \end{split}
\end{equation*}
where the true solution is 
\begin{equation*}
    u(t,x)=\frac{0.5}{\cosh(2(x-\sqrt{3}t))}-\frac{0.5}{\cosh(2(x-10+\sqrt{3}t))}+\frac{0.5}{\cosh(2(x+\sqrt{3}t))}-\frac{0.5}{\cosh(2(x+10-\sqrt{3}t))}.
\end{equation*}

\begin{figure}[htbp]
    \centering
    \includegraphics[height=0.2\textheight]{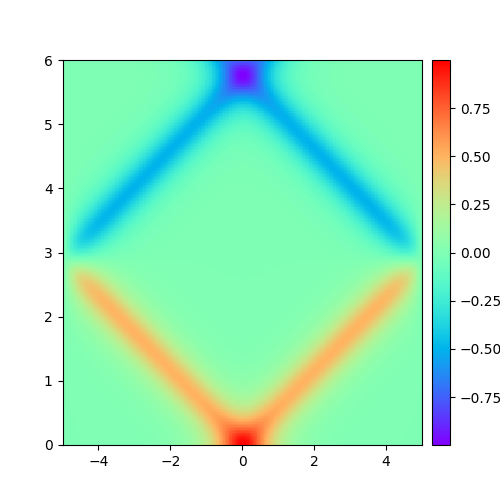}
    \caption{The plot of the exact solution.}
    \label{timedependant-true}
\end{figure}

\begin{figure}
	\centering
	\subcaptionbox*{R-PINN}
	{
		\includegraphics[height=0.2\textheight]{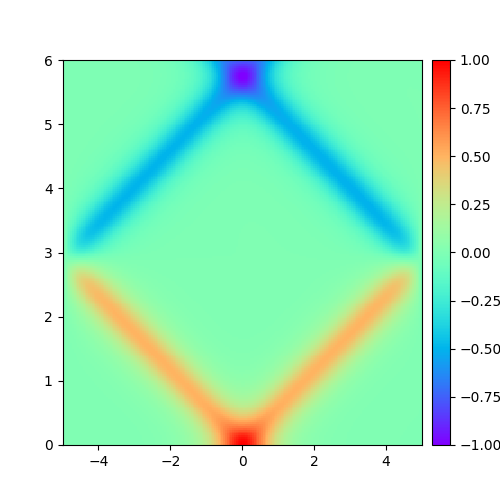}
        
        \includegraphics[height=0.2\textheight]{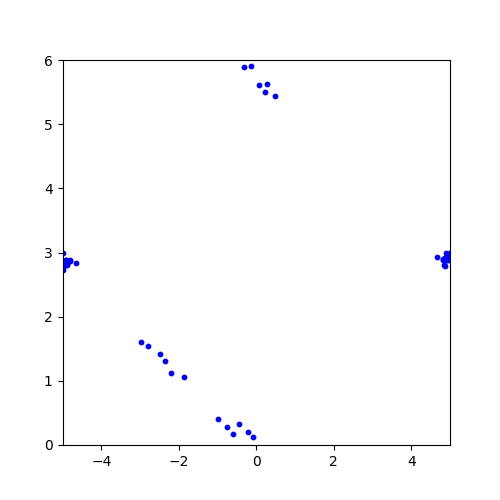}
	}
	\subcaptionbox*{FI-PINN}
	{
		\includegraphics[height=0.2\textheight]{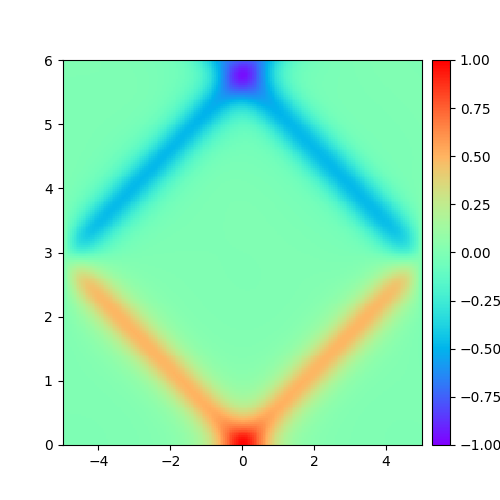}
        
		\includegraphics[height=0.2\textheight]{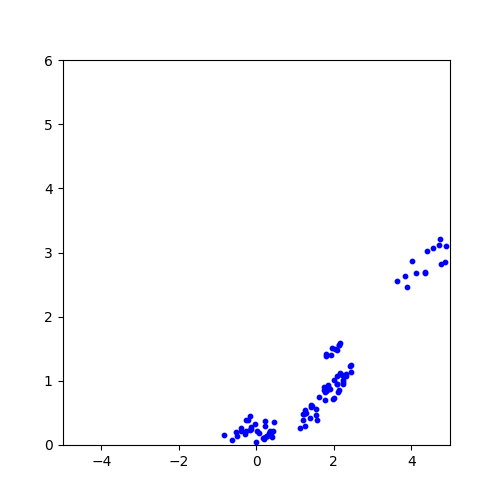}
	}
	\caption{The predicted solutions (left) and the distribution of adaptively distributed points (right) of two adaptive sampling methods solving the wave equation.}
	\label{timedependent-err-solution}
\end{figure}

\begin{figure}
    \centering
    \includegraphics[height=0.2\textheight]{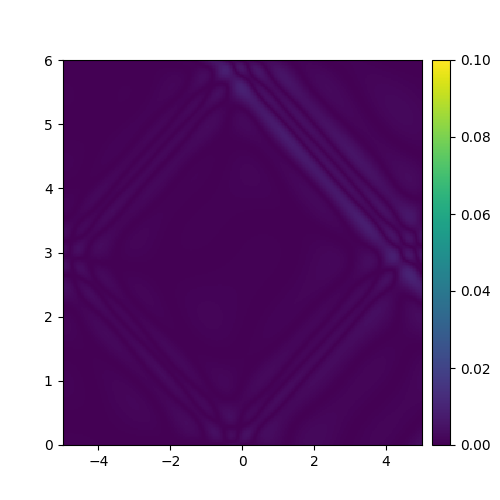}
    \includegraphics[height=0.2\textheight]{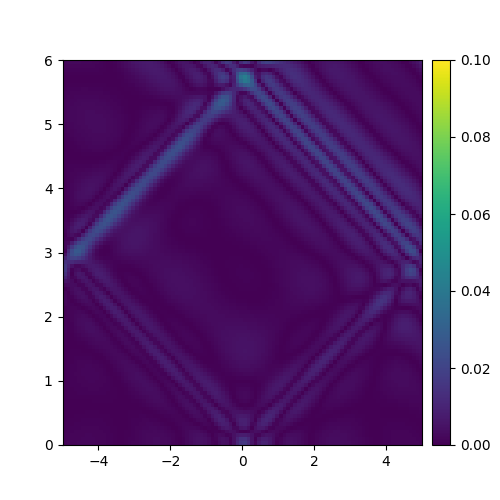}
    \includegraphics[height=0.2\textheight]{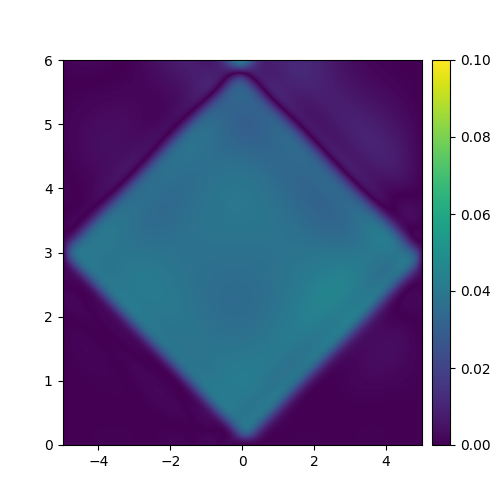}
    \caption{Comparison of absolute errors distributions across the domain for three different PINNs: R-PINN (left), FI-PINN (middle), and RBA-PINN (right).}
    \label{Timedependent-err}
\end{figure}

\begin{figure}[htbp]
    \centering
    \includegraphics[height=0.28\textheight]{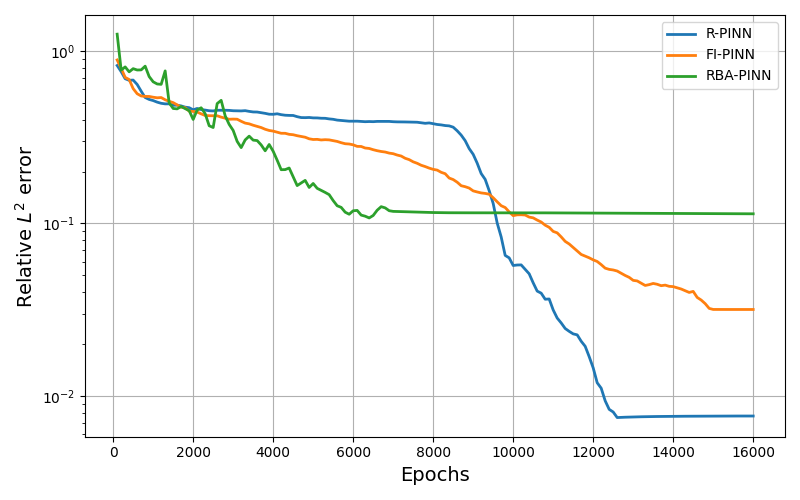}
    \caption{The relative $L^2$ error by three different PINNs in the optimization process of solving the wave equation.}
    \label{timedependent-L2}
\end{figure}

In this test, we adopt the adaptive iteration process with the number of $M=4$ and pre-train the network for 5000 iterations. We show the relative $L^2$ error and $L^\infty$ error in the Table~\ref{timedependent-l2}.
The relative $L^2$ error and $L^\infty$ error can reach 0.0076 and 0.0098 with R-PINN, which are obviously lower than the other two methods. 
The decreasing trend of the relative $L^2$ error during training can be clearly observed in Fig.~\ref{timedependent-L2}. The exact solution is depicted in Fig.~\ref{timedependant-true}. The left panel of Fig.~\ref{timedependent-err-solution} compares the predicted solutions of the two adaptive sampling methods. Fig.~\ref{Timedependent-err} illustrates the absolute errors across the domain, showing that the errors associated with R-PINN are smaller than those with FI-PINN and  RBA-PINN. Finally, the right panel of Fig.~\ref{timedependent-err-solution} displays the adaptively distributed points for two adaptive sampling methods. From the right panel, it can be observed that FI-PINN encounters similar difficulties when resolving multiple regions with large errors, indicating its limited effectiveness in such scenarios, where R-PINN demonstrates better performance.


\begin{figure}[b]
    \centering
    \includegraphics[height=0.24\textheight]{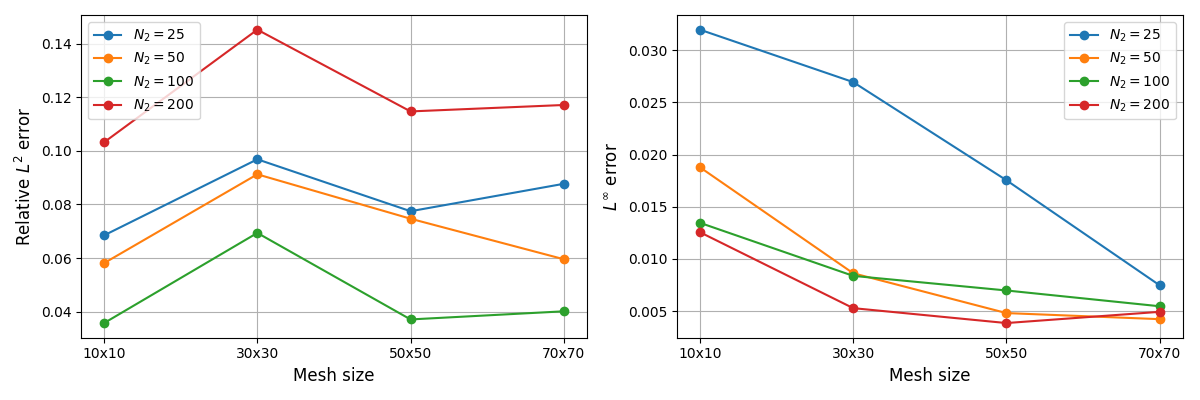}
    \caption{Comparison of the relative $L^2$ error (left) and $L^\infty$
  error (right) under varying mesh sizes for different values of $N_2$.}
    \label{mesh-err}
\end{figure}

\subsection{Tests on the effects of sampling and mesh}

For R-PINN, two critical factors influencing accuracy are the number of collocation points and the mesh size. To evaluate their effects, we conduct experiments on the Poisson's equation presented in Section~\ref{poisson} and examine how these two aspects affect the performance of R-PINN.

The collocation points consist of two components: the background collocation points and the adaptively distributed points. We begin by analyzing the combined influence of $N_2$ and the mesh size. Specifically, we fix $N_1 = 8000$ and vary the mesh size among $10\times10$, $30\times30$, $50\times50$, and $70\times70$, while also testing $N_2 = 25,\ 50,\ 100,\ 200$. The finest mesh used is $70\times70$, which already generates 9800 elements. A finer mesh would result in only a small number of elements containing background collocation points, thus undermining the effectiveness of the method.

We first examine the effect of mesh size, illustrated in Fig.~\ref{mesh}. A decreasing trend in the $L^\infty$ error is clearly observed with finer meshes, whereas the relative $L^2$ error does not follow a consistent trend. This is reasonable because the mesh is primarily used to guide a-posteriori error estimation and sampling, rather than to enrich the representational capacity of the solution space.
Next, we analyze the impact of $N_2$, which can also be seen in Fig.~\ref{mesh-err} by comparing four different curves in each subfigure. Both the relative $L^2$ error and the $L^\infty$  error generally decrease as $N_2$ increases, with the exception of a significant increase in relative $L^2$ error when $N_2$ increases from 100 to 200. This anomaly results from an imbalance between the background collocation points and adaptively distributed points, where excessive adaptively distributed points destroy the intended sampling strategy.

\begin{figure}
    \centering
    \includegraphics[height=0.24\textheight]{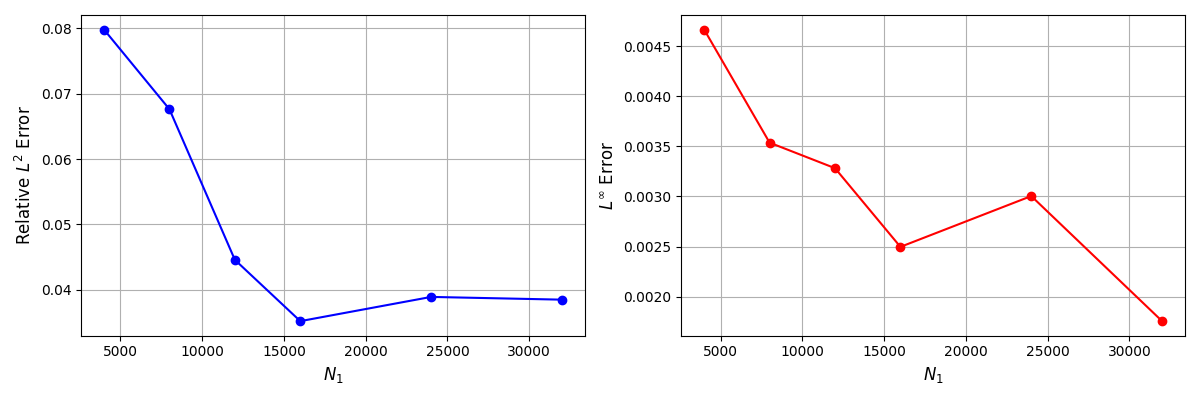}
    \caption{Comparison of the relative $L^2$ error (left) and $L^\infty$ error (right) under varying values of $N_1$.}
    \label{back-eff}
\end{figure}

Finally, we examine the effect of the number of background collocation points $N_1$. With the mesh size fixed at $50\times50$ and $N_2 = 100$, we evaluate the error under different values of $N_1$, as shown in Fig.~\ref{back-eff}. It is clearly observed that both the relative $L^2$ and $L^\infty$ errors decrease with increasing $N_1$, indicating that more background points contribute to a more accurate solution.

\section{Conclusion}

In this paper, we present recovery-type a-posteriori estimator enhanced adaptive PINN that combines the PINN with a recovery type estimator which is popularly used in the adaptive finite element method. Our innovation lies in the use of the recovery-type a-posteriori error estimator, which is different from the residual of PINNs, and the design of a tailored sampling method RecAD. We compare the performance of the proposed R-PINN with FI-PINN and RBA-PINN on several problems featuring sharp gradients or singularities, which vanilla PINNs struggle to resolve effectively. The results demonstrate that R-PINN achieves significant improvements. More studies can be done to analyze the performance of R-PINN for more complex problems and we shall focus on and analyze these issues in our future work.

\section{Acknowledgment}

JJ was partially supported by China Natural National Science Foundation (No. 22341302), the Fundamental Research Funds for the Central Universities, JLU (93Z172023Z05), and the Key Laboratory of Symbolic Computation and Knowledge Engineering of Ministry of Education of China housed at Jilin University. YL acknowledged the support of NSF-DMS (No. 2208499). CZ was partially supported by the Strategic Priority Research Program of the Chinese Academy of Sciences (No. XDB0640000) and the National Natural Science Foundation of China (No. 12571445). The authors wish to thank the anonymous reviewers for their constructive comments and valuable suggestions, which significantly improved the quality of this manuscript.

 \bibliographystyle{elsarticle-num-names} 
 \bibliography{paper}






\end{document}